\documentclass[10.5pt]{article}
\usepackage{amsfonts}
\usepackage{graphicx}

\usepackage{amsmath}
\usepackage{amssymb}
\usepackage{latexsym}
\usepackage{amsmath, amsfonts,amssymb, amsthm, euscript,makeidx,color,mathrsfs}

\def\sqr#1#2{{\vcenter{\vbox{\hrule height.#2pt
              \hbox{\vrule width.#2pt height#1pt \kern#1pt \vrule width.#2pt}
              \hrule height.#2pt}}}}
\def\signed #1{{\unskip\nobreak\hfil\penalty50
              \hskip2em\hbox{}\nobreak\hfil#1
              \parfillskip=0pt \finalhyphendemerits=0 \par}}
\def\endpf{\signed {$\sqr69$}}

\def\3n{\negthinspace \negthinspace \negthinspace }
\def\2n{\negthinspace \negthinspace }
\def\1n{\negthinspace }

\def\dbE{\mathbb{E}}
\def\dbF{\mathbb{F}}

\def\dbH{\mathbb{H}}

\def\dbP{\mathbb{P}}

\def\dbR{\mathbb{R}}
\def\dbS{\mathbb{S}}

\def\sL{\mathscr{L}}
\def\sM{\mathscr{M}}
\def\sN{\mathscr{N}}

\def\sP{\mathscr{P}}

\def\sS{\mathscr{S}}

\def\={\buildrel \triangle \over =}

\def\e{\varepsilon}
\def\z{\zeta}

\def\n{\nu}
\def\si{\sigma}

\def\f{\varphi}

\def\G{\Gamma}
\def\D{\Delta}
\def\Th{\Theta}
\def\L{\Lambda}

\def\cF{{\cal F}}

\def\cR{{\cal R}}

\def\cU{{\cal U}}

\def\cX{{\cal X}}

\def\no{\noindent}

\def\ms{\medskip}

\def\q{\quad}
\def\qq{\qquad}
\def\hb{\hbox}

\def\lan{\mathop{\langle}}
\def\ran{\mathop{\rangle}}

\def\wt{\widetilde}

\def\cd{\cdot}

\def\ae{\hbox{\rm a.e.{ }}}
\def\as{\hbox{\rm a.s.{ }}}

\def\deq{\mathop{\buildrel\D\over=}}
\def\les{\leqslant}
\def\ges{\geqslant}

\def\({\Big (}
\def\){\Big )}
\def\[{\Big[}
\def\]{\Big]}
\def\bde{\begin{definition}}
\def\ede{\end{definition}}
\def\be{\begin{equation}}
\def\bel{\begin{equation}\label}
\def\ee{\end{equation}}
\def\bt{\begin{theorem}}
\def\et{\end{theorem}}
\def\bc{\begin{corollary}}
\def\ec{\end{corollary}}
\def\bl{\begin{lemma}}
\def\el{\end{lemma}}
\def\bp{\begin{proposition}}
\def\ep{\end{proposition}}
\def\bas{\begin{assumption}}
\def\eas{\end{assumption}}
\def\br{\begin{remark}}
\def\er{\end{remark}}
\def\ba{\begin{array}}
\def\ea{\end{array}}
\def\ed{\end{document}}

\def\square#1{\vbox{\hrule\hbox{\vrule height#1%
     \kern#1\vrule}\hrule}}
\def\rectangle#1#2{\vbox{\hrule\hbox{\vrule height#1%
     \kern#2\vrule}\hrule}}

\font\tenbb=msbm10 \font\sevenbb=msbm7 \font\fivebb=msbm5

\newfam\bbfam
\scriptscriptfont\bbfam=\fivebb \textfont\bbfam=\tenbb
\scriptfont\bbfam=\sevenbb

\newtheorem{lemma}{Lemma}[section]
\newtheorem{remark}{Remark}[section]

\newtheorem{theorem}{Theorem}[section]
\newtheorem{corollary}{Corollary}[section]

\newtheorem{definition}{Definition}[section]
\newtheorem{proposition}{Proposition}[section]
\newtheorem{assumption}{Assumption}[section]

\makeatletter
   
   \@addtoreset{equation}{section}
\makeatother

\begin{document}

\title{\bf Linear Quadratic Stochastic Two-Person Zero-Sum Differential Games in an Infinite Horizon\thanks{
This work is supported in part by NSF Grant DMS-1007514, China Scholarship Council, and NSFC Grant 973-
2007CB814901.}}
\author{Jingrui Sun\thanks{School of Mathematical Sciences, University of Science and Technology of China, Hefei, 230026, P. R. China;
Email: sjr@mail.ustc.edu.cn.}\ , \ Jiongmin Yong\thanks{Department of Mathematics, University of Central Florida, Orlando, FL 32816, USA;
Email: jiongmin.yong@ucf.edu.}\ , \ and \ Shuguang Zhang\thanks{Department of Statistics and Finance, University of Science and Technology of China, Hefei, 230026, P. R. China; Email: sgzhang@ustc.edu.cn.}
 }
\maketitle

\no\bf Abstract: \rm This paper is concerned with a linear quadratic
stochastic two-person zero-sum differential game with constant
coefficients in an infinite time horizon. Open-loop and closed-loop
saddle points are introduced. The existence of closed-loop saddle
points is characterized by the solvability of an algebraic Riccati
equation with a certain stabilizing condition. A crucial result
makes our approach work is the unique solvability of a class of
linear backward stochastic differential equations in an infinite
horizon.

\ms

\no\bf Keywords: \rm linear quadratic stochastic differential game, two-person, zero-sum, infinite horizon,
open-loop and closed-loop saddle points, algebraic Riccati equation, stabilizing solution.

\ms

\no\bf AMS Mathematics Subject Classification. \rm \rm 93E20, 91A23,
49N10, 49N70.

\section{Introduction}

Let $(\Omega,\cF,\dbF,\dbP)$ be a complete filtered probability space on which a one-dimensional standard Brownian motion
$W(\cd)$ is defined with $\dbF=\{\cF_t\}_{t\ges0}$ being its natural filtration augmented by all the $\dbP$-null sets
in $\cF$ \cite{Karatzas-Shreve 1991,Yong-Zhou 1999}.
Consider the following controlled linear stochastic differential equation (SDE, for short) on the infinite time horizon
$[0,\infty)$:
\bel{state}\left\{\2n\ba{ll}
\noalign{\smallskip}\displaystyle  dX(t)=\big[AX(t)+B_1u_1(t)+B_2u_2(t)+b(t)\big]dt\\
\noalign{\smallskip}\displaystyle \qq\qq\q+\big[CX(t)+D_1u_1(t)+D_2u_2(t)+\si(t)\big]dW(t),\q t\ges0, \\
\noalign{\smallskip}\displaystyle  X(0)= x,\ea\right.\ee
where $A,C\1n\in\1n\dbR^{n\times n}$ and $B_i\1n\in\1n\dbR^{n\times m_i}$ $(i=1,2)$ are given (deterministic) matrices; $b(\cd)$
and $\si(\cd)$ are $\dbR^n$-valued, $\dbF$-adapted, square integrable processes. In the above, $X(\cd)$, valued in
$\dbR^n$, is called the {\it state process} with {\it initial state} $x\in\dbR^n$; for $i=1,2$, $u_i(\cd)$, valued in
$\dbR^{m_i}$, is called the {\it control process} of Player $i$. Let $\dbH$ be a Euclidean space and $T>0$ , we introduce
the following:
$$\ba{ll}
\noalign{\smallskip}\displaystyle  L_\dbF^2(\dbH)=\Big\{\f:[0,\infty)\times\Omega\to\dbH\bigm|\f(\cd)\hb{ is
$\dbF$-adapted, }\dbE\int^\infty_0|\f(t)|^2dt<\infty\Big\},\\
\noalign{\smallskip}\displaystyle  \cX[0,T]=\Big\{X:[0,\infty)\times\Omega\to\dbR^n\bigm|X(\cd)\hb{ is $\dbF$-adapted, continuous, }
\dbE\(\sup_{0\les t\les T}|X(t)|^2\)<\infty\Big\}, \\
\noalign{\smallskip}\displaystyle \cX_{loc}[0,\infty)=\bigcap_{T>0}\cX[0,T],\qq
\cX[0,\infty)=\Big\{X(\cd)\in\cX_{loc}[0,\infty)\bigm|
\dbE\int_0^\infty|X(t)|^2dt<\infty\Big\}.\ea$$
By a standard argument using contraction mapping theorem, one can show that for any initial state $x\in\dbR^n$ and control pair $(u_1(\cd),u_2(\cd))\1n\in\1n L_\dbF^2(\dbR^{m_1})\1n\times\1n L_\dbF^2(\dbR^{m_2})$, state equation (\ref{state}) admits a unique solution
$X(\cd)\equiv X(\cd\,;x,u_1(\cd),u_2(\cd))\in\cX_{loc}[0,\infty)$.
Next, we introduce the following performance functional:
\bel{cost}\ba{ll}
\noalign{\smallskip}\displaystyle  J(x;u_1(\cd),u_2(\cd))\\
\noalign{\smallskip}\displaystyle  \deq\dbE\int_0^\infty
\Big[\lan\begin{pmatrix}Q&S_1^T&S_2^T\\
                        S_1&R_{11}&R_{12}\\
                        S_2&R_{21}&R_{22}\end{pmatrix}
                        \begin{pmatrix}X(t)\\ u_1(t)\\ u_2(t)
                        \end{pmatrix},\begin{pmatrix}X(t)\\ u_1(t)\\ u_2(t)
                        \end{pmatrix}\ran
+2\lan\begin{pmatrix}q(t)\\ \rho_1(t)\\ \rho_2(t)\end{pmatrix},
                                          \begin{pmatrix}X(t)\\ u_1(t)\\ u_2(t)
                                          \end{pmatrix}\ran\Big]dt,\ea\ee
where
$$Q\1n\in\dbS^n,\  S_i\1n\in\dbR^{m_i\times n},\ R_{ii}\1n\in\dbS^{m_i},\ R_{21}^T\1n=\1nR_{12}\1n\in\dbR^{m_1\times m_2},
\  q(\cd)\1n\in L_\dbF^2(\dbR^n),\  \rho_i(\cd)\1n\in L_\dbF^2(\dbR^{m_i});\q i=1,2.$$
In the above, $\dbS^k$ is the set of all $(k\times k)$ symmetric
matrices, and $M^T$  is the transpose of $M$ (a matrix or a vector);
$X(\cd)=X(\cd\,;x,u_1(\cd),u_2(\cd))$ on the right hand side of
(\ref{cost}) is the corresponding state process. Note that in
general, for $(x,u_1(\cd),u_2(\cd))\in\dbR^n\times
L_\dbF^2(\dbR^{m_1})\times L_\dbF^2(\dbR^{m_2})$, the solution
$X(\cd)\equiv X(\cd\,;x,u_1(\cd),u_2(\cd))$ of (\ref{state}) might
just be in $\cX_{loc}[0,\infty)$ and the above performance
functional $J(x;u_1(\cd),u_2(\cd))$ might not be defined. Therefore,
we introduce the following set:
$$\cU_{ad}(x)\deq\big\{(u_1(\cd),u_2(\cd))\in L_\dbF^2(\dbR^{m_1})\1n\times\1nL_\dbF^2(\dbR^{m_2})
\bigm|X(\cd\,;x,u_1(\cd),u_2(\cd))\in\cX[0,\infty)\big\},\q x\in\dbR^n.$$
Any element $(u_1(\cd),u_2(\cd))\1n\in\cU_{ad}(x)$ is called an {\it
admissible control pair} for the initial state $x$ and the
corresponding $X(\cd)=X(\cd\,;x,u_1(\cd),u_2(\cd))$ is called an
{\it admissible state process} with the initial state $x$. Roughly
speaking, in the game, Player 1 wishes to minimize (\ref{cost}) by
selecting a control $u_1(\cd)$, and Player 2 wishes to maximize
(\ref{cost}) by selecting a control $u_2(\cd)$. Therefore,
(\ref{cost}) represents the cost for Player 1 and the payoff for
Player 2. The problem is to find an admissible control pair
$(u_1^*(\cd),u_2^*(\cd))$ that both players can accept, and we refer
to such a problem as a linear quadratic (LQ, for short) stochastic
{\it two-person zero-sum differential game}, denoted by Problem
(LQG). There are basically two types of controls for both players:
open-loop controls and closed-loop controls. An open-loop control
usually depends on the initial state as well as all the information,
including those of the opponent, over the whole time duration
$[0,\infty)$, whereas a closed-loop control is required to be
independent of the initial state, and the future information. Thus,
in reality, it is more meaningful and convenient to using
closed-loop controls rather than open-loop controls. However,
mathematically, open-loop controls are still meaningful and they are
actually helpful in finding ``optimal'' closed-loop controls.

\ms

Let us briefly recall some relevant history. In 1965, deterministic
LQ two-person zero-sum differential games in finite horizon (LQDG
problem, for short) was introduced and studied by Ho--Bryson--Baron
\cite{Ho-Bryson-Baron 1965}. In 1970, Schmitendorf studied both
open-loop and closed-loop strategies for LQDG problems
(\cite{Schmitendorf 1970}). Among other things, it was shown that
the existence of a closed-loop saddle point may not imply that of an
open-loop saddle point. In 1979, Bernhard carefully investigated
LQDG problems from closed-loop point of view (\cite{Bernhard 1979});
see also the book by Basar and Bernhard \cite{Basar-Bernhard 1991}
in this aspect. In 2005, Zhang \cite{Zhang 2004} proved that for an
LQDG problem, the existence of the open-loop value is equivalent to
the finiteness of the corresponding open-loop lower and upper
values, which is also equivalent to the existence of an open-loop
saddle point. Along this line, there were a couple of follow-up
works \cite{Delfour 2007,Delfour-Sbarba 2009} appeared afterwards.
In 2006, Mou--Yong studied a stochastic LQ two-person zero-sum
differential game in finite horizon from an open-loop point of
view, by means of Hilbert space method (\cite{Mou-Yong 2006}). On
the other hand, in 1976, Ichikawa studied a deterministic LQ
two-person zero-sum differential games on $[0,\infty)$ in a Hilbert
space and deduced some sufficient conditions for the existence of
closed-loop saddle points (\cite{Ichikawa 1976}). In 2000, Ait
Rami--Moore--Zhou studied an LQ stochastic optimal control problem
on $[0,\infty)$ (\cite{Rami-Zhou-Moore 2000}), followed by the work
of Wu--Zhou (\cite{Wu-Zhou 2001}). Recently, based on the work of Yong
\cite{Yong 2013}, Huang--Li--Yong studied a mean--field LQ optimal
control problem on $[0,\infty)$ (\cite{Huang-Li-Yong}).

\ms

The rest of the paper is organized as follows. In Section 2, we
collect some preliminary results. Section 3 is devoted to the unique
solvability of a linear backward stochastic differential equation (BSDE, for short) on
$[0,\infty)$. In Section 4, we discuss closed-loop optimal controls
of Problem (LQ) and deduce a necessary condition for the existence
of a closed-loop optimal control via the solvability of an algebraic
Riccati equation (ARE, for short). In Section 5, we pose our differential game
problem and characterize closed-loop saddle points by means of
algebraic Riccati equations. Some examples are presented in Section 6.

\section{Preliminary Results}

Let us begin by considering a stochastic optimal control problem.
The state equation takes the following form:
\bel{H-state}\left\{\2n\ba{ll}
\noalign{\smallskip}\displaystyle  dX(t)=\big[AX(t)+Bu(t)+b(t)\big]dt+\big[CX(t)+Du(t)+\si(t)\big]dW(t),\q t\ges0, \\
\noalign{\smallskip}\displaystyle  X(0)= x,\ea\right.\ee
with cost functional
\bel{H-cost}
J(x;u(\cd))
=\dbE\int_0^\infty\Big[\lan\begin{pmatrix}Q&S^T\\
                                      S&R\end{pmatrix}\begin{pmatrix}
                                      X(t)\\ u(t)\end{pmatrix},
                                      \begin{pmatrix}
                                      X(t)\\
                                      u(t)\end{pmatrix}\ran+
                                      2\lan\begin{pmatrix}q(t)\\ \rho(t)
                                      \end{pmatrix},\begin{pmatrix}X(t)\\
                                      u(t)\
                                      \end{pmatrix}\ran\Big] dt,\ee
where $A,C\1n\in\1n\dbR^{n\times n}$, $B,D\1n\in\1n\dbR^{n\times
m}$, $Q\1n\in\1n\dbS^n$, $R\1n\in\1n\dbS^m$,
$S\1n\in\1n\dbR^{m\times n}$, and $b(\cd),\si(\cd),q(\cd)\1n\in\1n
L_\dbF^2(\dbR^n)$, $\rho(\cd)\1n\in\1n L_\dbF^2(\dbR^m)$. The
solution of (\ref{H-state}) is denoted by $X(\cd\,;x,u(\cd))$. For
any given $x\in\dbR^n$, the set of {\it admissible controls} is
defined by the following:
$$\cU_{ad}(x)\deq\Big\{u(\cd)\in L^2_\dbF(\dbR^m)\bigm|X(\cd\,;x,u(\cd))\in\cX[0,\infty)\Big\}.$$
Clearly, $\cU_{ad}(x)$ is a convex subset of $L^2_\dbF(\dbR^m)$, but not a subspace of $L^2_\dbF(\dbR^m)$ in general.
We pose the following problem.

\ms

\bf Problem (LQ). \rm For any $x\in\dbR^n$, find a $\bar
u(\cd)\in\cU_{ad}(x)$, such that
\bel{2.4}V(x)\deq J(x;\bar u(\cd))=\inf_{u(\cd)\in\cU_{ad}(x)}J(x;u(\cd)).\ee
Any $\bar u(\cd)\in\cU_{ad}(x)$ satisfying (\ref{2.4}) is called an
{\it open-loop optimal control} of Problem (LQ), and the
corresponding $\bar X(\cd)\equiv X(\cd\,;x,\bar u(\cd))$ is called
an {\it optimal state process}. The function $V(\cd)$ is called the {\it
value function} of Problem (LQ). The following notions are similar
to those introduced in \cite{Yong-Zhou 1999}.

\ms

\bf Definition 2.1. \rm (i) Problem (LQ) is said to be {\it finite}
if
\bel{}V(x)>-\infty,\qq\forall x\in\dbR^n.\ee

(ii) Problem (LQ) is said to be ({\it uniquely}) {\it solvable} if
it has a (unique) open-loop optimal control.

\ms

When $b(\cd),\si(\cd)\1n=\1n0$, we briefly denote the system
(\ref{H-state}) by $[A,C;B,D]$. We also denote by $[A,C]$ the
following uncontrolled system:
\bel{AC}\left\{\2n\ba{ll}
\noalign{\smallskip}\displaystyle  dX(t)=AX(t)dt+CX(t)dW(t),\qq t\ges0, \\
\noalign{\smallskip}\displaystyle  X(0)= x.\ea\right.\ee
When $b(\cd),\si(\cd),q(\cd),\rho(\cd)\1n=\1n0$, we denote the
corresponding Problem (LQ) by Problem $\hb{(LQ)}^0$. The
corresponding cost functional and value function are denoted by
$J^0(x;u(\cd))$ and $V^0(x)$, respectively.

\ms

We note that, in general, the admissible control set $\cU_{ad}(x)$
may be empty for some $x\in\dbR^n$. To avoid such a case, we
introduce the following definition.

\ms

\bf Definition 2.2. \rm (i) System $[A,C]$ is said to be $L^2$-{\it
exponentially stable} if for any $x\in\dbR^n$, the solution
$X(\cd)\equiv X(\cd\,;x)\in\cX_{loc}[0,\infty)$ of (\ref{AC})
satisfies the following:
$$\lim_{t\to\infty}e^{\lambda t}\dbE|X(t)|^2=0,\q \hb{for some } \lambda>0.$$

\ms

(ii) System $[A,C]$ is said to be $L^2$-{\it globally integrable} if for any $x\in\dbR^n$, the solution
$X(\cd)\equiv X(\cd\,;x)\in\cX_{loc}[0,\infty)$ of (\ref{AC}) is in $\cX[0,\infty)$.

\ms

(iii) System $[A,C]$ is said to be $L^2$-{\it asymptotically stable} if for any $x\in\dbR^n$, the solution
$X(\cd)\equiv X(\cd\,;x)\in\cX_{loc}[0,\infty)$ of (\ref{AC}) satisfies the following:
$$\lim_{t\to\infty}\dbE|X(t)|^2=0.$$

\ms

The following result will be used frequently in this paper. For a
proof, see \cite{Huang-Li-Yong}.

\ms

\bf Lemma 2.3. \sl The following are equivalent:

\ms

{\rm(i)} System $[A,C]$ is $L^2$-exponentially stable;

\ms

{\rm(ii)} System $[A,C]$ is $L^2$-globally integrable;

\ms

{\rm(iii)} For any $\Lambda>0$, the following Lyapunov equation admits a
solution $P>0$:
\bel{Lyapunov}PA+A^TP+C^TPC+\Lambda=0;\ee

\ms

{\rm(iv)} There exists a $P>0$ such that $PA+A^TP+C^TPC<0$;

\ms

{\rm(v)} System $[A,C]$ is $L^2$-asymptotically stable, and there exists a $P\in\dbS^n$ such that
$$PA+A^TP+C^TPC<0.$$
In this case, we simply say that the system $[A,C]$ is $L^2$-stable.

\ms

\rm

Next, we present a result concerning the $L^2$-integrability of the
solution to the following system:
\bel{L2-int}\left\{\2n\ba{ll}
\noalign{\smallskip}\displaystyle  dX(t)=\big[AX(t)+b(t)\big]dt+\big[CX(t)+\si(t)\big]dW(t),\q t\ges0, \\
\noalign{\smallskip}\displaystyle  X(0)= x.\ea\right.\ee

\ms

\bf Proposition 2.4. \sl Let $A,C\in\dbR^{n\times n}$ and
$b(\cd),\si(\cd)\1n\in\1n L^2_\dbF(\dbR^n)$. Let $X(\cd)\equiv
X(\cd\,;x)$ be the solution to the SDE $(\ref{L2-int})$. If $[A,C]$
is $L^2$-stable, then $X(\cd)\in\cX[0,\infty)$.

\ms

\it Proof. \rm Since $[A,C]$ is $L^2$-stable, by Lemma 2.3, there
exists a $P>0$ such that
$$PA+A^TP+C^TPC\equiv-\Lambda<0.$$
Applying It\^o's formula to $s\mapsto \lan PX(s),X(s)\ran$, one has
$$\ba{ll}
\noalign{\smallskip}\displaystyle  \dbE\lan PX(t),X(t)\ran-\lan Px,x\ran\\
\noalign{\smallskip}\displaystyle \,=\dbE\int_0^t\1n\Big[\lan\big(PA\1n+\1nA^TP\1n+\1nC^TPC\big)X(s),X(s)\ran\\
\noalign{\smallskip}\displaystyle \qq\qq\q+2\lan Pb(s)\1n+\1nC^TP\si(s),X(s)\ran\1n+\1n\lan P\si(s),\si(s)\ran\Big]ds\\
\noalign{\smallskip}\displaystyle \,=\dbE\int_0^t\1n\Big[-\1n\lan\Lambda X(s),X(s)\ran\1n+2\lan
Pb(s)\1n +\1nC^TP\si(s),X(s)\ran\1n+\1n\lan
P\si(s),\si(s)\ran\Big]ds,\q \forall t\ges0.\ea$$
Therefore
$$\ba{ll}
\noalign{\smallskip}\displaystyle  {d\over dt}\dbE\lan P^{1\over2}X(t),P^{1\over2}X(t)\ran={d\over dt}\dbE\lan PX(t),X(t)\ran\\
\noalign{\smallskip}\displaystyle \,=-\dbE\lan\Lambda X(t),X(t)\ran+2\dbE\lan Pb(t)\1n+\1nC^TP\si(t),X(t)\ran+\dbE\lan P\si(t),\si(t)\ran\\
\noalign{\smallskip}\displaystyle \,=-\dbE\lan\G
P^{1\over2}X(t),P^{1\over2}X(t)\ran+2\dbE\lan\eta(t),P^{1\over2}X(t)\ran+\dbE\lan
P\si(t),\si(t)\ran,\ea$$
where
$$\G\deq P^{-{1\over2}}\Lambda P^{-{1\over2}}>0,
\q \eta(\cd)=P^{1\over2}b(\cd)+P^{-{1\over2}}C^TP\si(\cd).$$
Let $\lambda>0$ be the smallest eigenvalue of $\G$. By Cauchy--Schwarz's
inequality, we have
$$\ba{ll}
\noalign{\smallskip}\displaystyle  {d\over dt}\dbE\lan P^{1\over2}X(t),P^{1\over2}X(t)\ran\\
\noalign{\smallskip}\displaystyle \,\les-\lambda\dbE\lan P^{1\over2}X(t),P^{1\over2}X(t)\ran+{\lambda\over2}\dbE\lan P^{1\over2}X(t),P^{1\over2}X(t)\ran
+{2\over\lambda}\dbE|\eta(t)|^2+\dbE\lan P\si(t),\si(t)\ran\\
\noalign{\smallskip}\displaystyle \,=-{\lambda\over2}\dbE\lan P^{1\over2}X(t),P^{1\over2}X(t)\ran+{2\over\lambda}\dbE|\eta(t)|^2+\dbE\lan P\si(t),\si(t)\ran.\ea$$
Let $\mu>0$ be the smallest eigenvalue of $P$. By Gronwall's
inequality, we obtain
$$\ba{ll}
\noalign{\smallskip}\displaystyle  \mu\dbE|X(t)|^2\les\dbE\lan P^{1\over2}X(t),P^{1\over2}X(t)\ran\\
\noalign{\smallskip}\displaystyle \qq\qq\q\2n\les e^{-{\lambda\over2}t}\lan Px,x\ran
+\int_0^te^{-{\lambda\over2}(t-s)}\big[\,{2\over\lambda}\dbE|\eta(s)|^2+\dbE\lan P\si(s),\si(s)\ran\big]ds,\ea$$
which, together with Young's inequality, implies that $\dbE|X(\cd)|^2$ is integrable over $[0,\infty)$. \endpf

\ms

\bf Definition 2.5. \rm System $[A,C;B,D]$ is said to be {\it
$L^2$-stabilizable} if there exists a $\Th\1n\in\1n\dbR^{m\times n}$
such that $[A+B\Th,C+D\Th]$ is $L^2$-stable. In this case, $\Th$ is
called a {\it stabilizer} of $[A,C;B,D]$. We denote the set of all
stabilizers of $[A,C;B,D]$ by $\sS\equiv\sS[A,C;B,D]$.

\ms

We now introduce the following assumption.

\ms

{\bf(H1)} System $[A,C;B,D]$ is $L^2$-stabilizable, i.e.,
\bel{}\sS[A,C;B,D]\ne\varnothing.\ee

\ms

By Proposition 2.4, we see that under (H1), $\cU_{ad}(x)$ is
nonempty for any $x\in\dbR^n$. Moreover, we have the following
proposition.

\ms

\bf Proposition 2.6. \sl Let {\rm(H1)} hold. Then for any
$x\in\dbR^n$, $u(\cd)\1n\in\1n\cU_{ad}(x)$ if and only if
\bel{}u(\cd)=\Th X(\cd)+v(\cd),\ee
for some $\Th\in\sS[A,C;B,D]$ and $v(\cd)\in L^2_\dbF(\dbR^m)$,
where $X(\cd)$ is the solution of the following SDE:
\bel{A2-1}\left\{\2n\ba{ll}
\noalign{\smallskip}\displaystyle  dX(t)=\big[(A+B\Th)X(t)+Bv(t)+b(t)\big]dt\\
\noalign{\smallskip}\displaystyle \qq\qq\q+\big[(C+D\Th)X(t)+Dv(t)+\si(t)\big]dW(t),\q t\ges0, \\
\noalign{\smallskip}\displaystyle  X(0)= x.\ea\right.\ee

\ms

\it Proof. \rm Let $v(\cd)\in L^2_\dbF(\dbR^m)$ and $X(\cd)$ be the
solution of (\ref{A2-1}). Since $[A+B\Th,C+D\Th]$ is $L^2$-stable,
by Proposition 2.4, $X(\cd)\in\cX[0,\infty)$. Set
$$u(\cd)\deq\Th X(\cd)+v(\cd)\in L^2_\dbF(\dbR^m).$$
By uniqueness, $X(\cd)$ also solves the following SDE:
\bel{A2-2}\left\{\2n\ba{ll}
\noalign{\smallskip}\displaystyle
dX(t)=\big[AX(t)+Bu(t)+b(t)\big]dt+\big[CX(t)+Du(t)+\si(t)\big]dW(t),\q t\ges 0,\\
\noalign{\smallskip}\displaystyle  X(0)=x.\ea\right.\ee
Thus, $u(\cd)\in\cU_{ad}(x)$.

\ms

On the other hand, suppose $u(\cd)\in\cU_{ad}(x)$. Let
$X(\cd)\in\cX[0,\infty)$ be the solution of (\ref{A2-2}). Pick any
$\Th\in\sS[A,C;B,D]$ and set
$$v(\cd)\deq u(\cd)-\Th X(\cd)\in L^2_\dbF(\dbR^m).$$
By uniqueness, $X(\cd)$ also solves (\ref{A2-1}). Thus, $u(\cd)=\Th
X(\cd)+v(\cd)$ with $X(\cd)$ being the solution of (\ref{A2-1}).
\endpf

\ms

Now, we introduce the following notations:
$$\sM(P)\1n=\1nPA\1n+\1nA^TP\1n+\1nC^TPC\1n+\1nQ,~~\sL(P)\1n=\1nPB\1n+\1nC^TPD\1n+\1nS^T\2n,~~
\sN(P)\1n=\1nR\1n+\1nD^TPD,\q \forall P\in\dbS^n,$$
and define the following convex set:
$$\sP\deq\bigg\{P\in\dbS^n\biggm|\begin{pmatrix}\sM(P)&\sL(P)\\\sL(P)^T&\sN(P)\end{pmatrix}\ges0\bigg\}.$$

\ms

The following result, found in \cite{Rami-Zhou-Moore 2000},
characterizes the finiteness of Problem $\hb{(LQ)}^0$.

\ms

\bf Lemma 2.7. \sl Problem $\hb{\rm(LQ)}^0$ is finite if and only if
$\sP\neq\varnothing$. In this case, $\sP$ has a maximal element
$P\in\sP$ (i.e., $P\ges\wt P$ $\forall \wt P\in\sP$). Moreover, we
have
$$V^0(x)=\lan Px,x\ran,\qq\forall x\in\dbR^n.$$

\rm

\section{Linear BSDEs in an Infinite Horizon}

In this section, we consider the following BSDE in the infinite time horizon $[0,\infty)$:
\bel{BSDE}dY(t)=-\big[A^TY(t)+C^TZ(t)+\f(t)\big]dt+Z(t)dW(t),\q t\in[0,\infty).\ee

\ms

\bf Definition 3.1. \rm An $L^2$-{\it stable adapted solution} of
(\ref{BSDE}) is a pair
$(Y(\cd),Z(\cd))\1n\in\cX[0,\infty)\1n\times\1n L_\dbF^2(\dbR^n)$
satisfying
\bel{3.2}Y(t)=Y(0)\1n-\1n\int_0^t\1n\big[A^TY(s)\1n+\1nC^TZ(s\1n)+\1n\f(s)\big]ds\1n+\1n\int_0^t\1nZ(s)dW(s),\q
\forall t\in[0,\infty),\q\as\ee

\ms

Note that by (\ref{3.2}), for any $T\in[0,\infty)$,
\bel{}Y(t)=Y(T)\1n+\1n\int_t^T\1n\big[A^TY(s)\1n+\1nC^TZ(s)\1n+\1n\f(s)\big]ds\1n-\1n\int_t^T\1nZ(s)dW(s),\qq
t\in[0,T],\q\as\ee
Hence, letting $T\to\infty$, we have
\bel{}Y(t)=\int_t^\infty\1n\big[A^TY(s)\1n+\1nC^TZ(s)\1n+\1n\f(s)\big]ds\1n-\1n\int_t^\infty\1nZ(s)dW(s),\qq t\in[0,\infty),\q\as\ee
This is a familiar form of linear BSDE on $[0,\infty)$. In 2000,
Peng and Shi considered the following BSDE:
\bel{BSDE1}dY(t)=-\big[G\big(t,Y(t),Z(t)\big)+\f(t)\big]dt+Z(t)dW(t),\q t\in[0,\infty),\ee
and it was shown that, under some mild conditions, equation
(\ref{BSDE1}) admits a unique adapted solution $(Y(\cd),Z(\cd))$
(\cite[Theorem 4]{Peng-Shi 2000}). In terms of $L^2$-stable adapted
solutions of (\ref{BSDE}), we can restate the result of
\cite{Peng-Shi 2000} as follows.

\ms

\bf Proposition 3.2. \sl Suppose
\bel{3.6}A+A^T+C^TC<0.\ee
Then for any $\f(\cd)\1n\in\1n L^2_\dbF(\dbR^n)$, BSDE
$(\ref{BSDE})$ admits a unique $L^2$-stable adapted solution
$(Y(\cd),Z(\cd))$.

\ms

\rm

Instead of the above, we have the following result which gives the
unique solvability of BSDE (\ref{BSDE1}) under a weaker condition.

\ms

\bf Theorem 3.3. \sl Suppose that $[A,C]$ is $L^2$-stable. Then for
any $\f(\cd)\1n\in\1n L^2_\dbF(\dbR^n)$, BSDE $(\ref{BSDE})$ admits
a unique $L^2$-stable adapted solution $(Y(\cd),Z(\cd))$.

\ms

\rm

Before proving the above result, let us make an observation. By
Lemma 2.3, part (iv), taking $P=I$, we see that condition
(\ref{3.6}) implies the $L^2$-stability of $[A,C]$. On the other
hand, let
$$A=\begin{pmatrix}-1&1\\-1&0\end{pmatrix},\q
C=\begin{pmatrix}{\sqrt{2}\over2}&0\\0&{\sqrt{2}\over2}\end{pmatrix},\q
P=\begin{pmatrix}2&-1\\-1&2\end{pmatrix}>0.$$
One has
$$PA+A^TP+C^TPC=\begin{pmatrix}-1&{1\over2}\\{1\over2}&-1\end{pmatrix}<0.$$
By Lemma 2.3, part (iv), $[A,C]$ is $L^2$-stable. However,
$$A+A^T+C^TC=\begin{pmatrix}-{3\over2}&0\\0&{1\over2}\end{pmatrix},$$
which is indefinite. Thus, (\ref{3.6}) fails. Hence, the condition
assumed in Theorem 3.3 is weaker than that assumed in Proposition
3.2. In order to prove Theorem 3.3, we need the following a priori
estimates.

\ms

\bf Proposition 3.4. \sl Suppose that $[A,C]$ is $L^2$-stable and
$\f(\cd)\1n\in\1n L^2_\dbF(\dbR^n)$. Let $(Y(\cd),Z(\cd))$ be an
$L^2$-stable adapted solution of BSDE $(\ref{BSDE})$. Then
\bel{E}\dbE\(\sup_{0\les
t<\infty}|Y(t)|^2\)+\dbE\int_0^\infty\1n|Z(t)|^2dt
\les K\dbE\int_0^\infty\1n|\f(t)|^2dt.\ee
Hereafter, $K>0$ represents a generic constant which can be
different from line to line.

\ms

\it Proof. \rm Since $[A,C]$ is $L^2$-stable,  by Lemma 2.3, there
exists a $P\1n>\1n0$ such that $PA\1n+\1nA^TP\1n+\1nC^TPC\1n<\1n0$.
Hence, one can choose $\e>0$ such that
$$PA+A^TP+(1+\e)C^TPC\equiv-\Lambda_\e<0.$$
Applying It\^o's formula to $s\1n\mapsto\1n \lan\1n
P^{-1}Y(s),Y(s)\1n\ran$, one has that for any $0\1n\les\1n
t\1n<\1nT\1n<\1n\infty$ (suppressing $s$ in the functions),
$$\ba{ll}
\noalign{\smallskip}\displaystyle  \lan P^{-1}Y(T),Y(T)\ran\1n-\1n\lan P^{-1}Y(t),Y(t)\ran\\
\noalign{\smallskip}\displaystyle \,=-\1n\int_t^T\2n\Big\{2\lan P^{-1}\big(A^TY\1n+\1nC^TZ\1n+\1n\f\big),Y\ran\1n-\1n\lan P^{-1}Z,Z\ran\1n\Big\}ds\1n
+\1n2\1n\int_t^T\2n\lan Z,P^{-1}Y\ran dW(s)\\
\noalign{\smallskip}\displaystyle \,=-\1n\int_t^T\2n\Big\{\1n\lan PAP^{-1}Y,P^{-1}Y\ran\1n+\1n\lan A^TPP^{-1}Y,P^{-1}Y\ran\1n+2\lan C^TZ,P^{-1}Y\ran\\
\noalign{\smallskip}\displaystyle \qq\qq\qq+2\lan\f,P^{-1}Y\ran\1n-\1n\lan P^{-1}Z,Z\ran\1n\Big\}ds\1n+\1n2\1n\int_t^T\2n\lan Z,P^{-1}Y\ran dW(s)\\
\noalign{\smallskip}\displaystyle \,=-\1n\int_t^T\2n\Big\{\1n\lan\big(PA+A^TP\big)P^{-1}Y,P^{-1}Y\ran\1n+2\lan\f,P^{-1}Y\ran\\
\noalign{\smallskip}\displaystyle \qq\qq\qq+2\lan Z,CP^{-1}Y\ran\1n-\1n\lan P^{-1}Z,Z\ran\1n\Big\}ds\1n+\1n2\1n\int_t^T\2n\lan Z,P^{-1}Y\ran dW(s)\\
\noalign{\smallskip}\displaystyle \,=-\1n\int_t^T\2n\Big\{\1n\lan-\Lambda_\e P^{-1}Y,P^{-1}Y\ran\1n
+2\lan\f,P^{-1}Y\ran\1n-(1+\e)\lan PCP^{-1}Y,CP^{-1}Y\ran\\
\noalign{\smallskip}\displaystyle \qq\qq\qq+2\lan Z,CP^{-1}Y\ran\1n-\1n\lan P^{-1}Z,Z\ran\1n\Big\}ds\1n+\1n2\1n\int_t^T\2n\lan Z,P^{-1}Y\ran dW(s)\\
\noalign{\smallskip}\displaystyle \,=-\1n\int_t^T\2n\Big\{\1n\lan-\Lambda_\e
P^{-1}Y,P^{-1}Y\ran\1n+2\lan\f,
P^{-1}Y\ran\\
\noalign{\smallskip}\displaystyle \qq\qq\qq-(1+\e)\lan P\big[CP^{-1}Y\1n
-{1\over1+\e}P^{-1}Z\big],CP^{-1}Y\1n-{1\over1+\e}P^{-1}Z\ran\\
\noalign{\smallskip}\displaystyle \qq\qq\qq-{\e\over1+\e}\lan
P^{-1}Z,Z\ran\1n\Big\}ds\1n+\1n2\1n\int_t^T\2n\lan Z,P^{-1}Y\ran
dW(s).\ea$$
Let $\lambda>0$ be the smallest eigenvalue of $\Lambda_\e\1n>0$. By
Cauchy--Schwarz's inequality, we have
\bel{E1}\ba{ll}
\noalign{\smallskip}\displaystyle \lan P^{-1}Y(t),Y(t)\ran\1n-\1n\lan P^{-1}Y(T),Y(T)\ran\1n
+\1n\int_t^T\2n{\e\over1+\e}\lan P^{-1}Z(s),Z(s)\ran ds\\
\noalign{\smallskip}\displaystyle \,=\int_t^T\2n\Big\{\1n\lan-\Lambda_\e P^{-1}Y,P^{-1}Y\ran\1n+2\lan\f,P^{-1}Y\ran\\
\noalign{\smallskip}\displaystyle \qq\qq-(1+\e)\big|P^{1\over2}\big[CP^{-1}Y\1n
-\1n{1\over1+\e}P^{-1}Z\big]\big|^2\Big\}ds
\1n-\1n2\1n\int_t^T\2n\lan Z,P^{-1}Y\ran dW(s)\\
\noalign{\smallskip}\displaystyle \,\les\int_t^T\2n\Big\{\1n-\1n\lambda|P^{-1}Y(s)|^2\1n
+\1n\lambda|P^{-1}Y(s)|^2\1n+\1n{1\over\lambda}|\f(s)|^2\Big\}ds
\1n-\1n2\1n\int_t^T\2n\lan Z(s),P^{-1}Y(s)\ran dW(s)\\
\noalign{\smallskip}\displaystyle \,={1\over\lambda}\int_t^T\1n|\f(s)|^2ds\1n-\1n2\1n\int_t^T\2n\lan
Z(s),P^{-1}Y(s)\ran dW(s).\ea\ee
Since $Y(\cd)\1n\in\1n\cX[0,\infty)$, we must have
$\lim_{T\to\infty}\dbE|Y(T)|^2=0$. Taking expectation on both sides
of (\ref{E1}), and letting $T\1n\to\1n\infty$, one has (noting that
$P>0$)
\bel{E2}\dbE|Y(t)|^2+\dbE\int_t^\infty|Z(s)|^2ds\les
K\dbE\int_t^\infty\1n|\f(s)|^2ds,\qq\forall t\in[0,\infty).\ee
On the other hand, by Burkholder--Davis--Gundy's inequality, we have
(noting (\ref{E2}))
\bel{E3}\ba{ll}
\noalign{\smallskip}\displaystyle \dbE\Big\{\sup_{0\les t\les T}\Big|\int_t^T\1n\lan Z(s),P^{-1}Y(s)\ran dW(s)
\Big|\,\Big\}\les2\dbE\Big\{\sup_{0\les t\les T}\Big|\int_0^t\lan\1n Z(s),P^{-1}Y(s)\ran dW(s)\Big|\,\Big\}\\
\noalign{\smallskip}\displaystyle \,\les K\dbE\Big\{\int_0^T\big|\lan Z(s),P^{-1}Y(s)\ran\big|^2ds\Big\}^{{1\over2}}
\les K\dbE\Big\{\int_0^T\big|P^{-{1\over2}}Z(s)\big|^2\big|P^{-{1\over2}}Y(s)\big|^2ds\Big\}^{{1\over2}}\\
\noalign{\smallskip}\displaystyle \,\les K\dbE\Big\{\Big(\sup_{0\les t\les T}\big|P^{-{1\over2}}Y(t)\big|^2\Big)^{{1\over2}}
\Big(\int_0^T\big|P^{-{1\over2}}Z(s)\big|^2ds\Big)^{{1\over2}}\Big\}\\
\noalign{\smallskip}\displaystyle \,\les {1\over4}\dbE\Big(\sup_{0\les t\les T}\big|P^{-{1\over2}}Y(t)\big|^2\Big)
+K\dbE\int_0^T\big|Z(s)\big|^2ds\\
\noalign{\smallskip}\displaystyle \,\les {1\over4}\dbE\Big(\sup_{0\les t\les
T}\big|P^{-{1\over2}}Y(t)\big|^2\Big)
+K\dbE\int_0^\infty\1n|\f(s)|^2ds.\ea\ee
Consequently, from (\ref{E1}), we obtain (using (\ref{E2})--(\ref{E3}))
$$\ba{ll}
\noalign{\smallskip}\displaystyle \dbE\Big(\sup_{0\les t\les T}\big|P^{-{1\over2}}Y(t)\big|^2\Big)
=\dbE\Big(\sup_{0\les t\les T}\lan P^{-1}Y(t),Y(t)\ran\Big)\\
\noalign{\smallskip}\displaystyle \,\les\dbE\lan\1n
P^{-1}Y(T),Y(T)\1n\ran\1n+{1\over\lambda}\dbE\int_0^T\1n|\f(s)|^2ds\1n
+\1n2\dbE\Big\{\sup_{0\les t\les T}\Big|\int_t^T\2n\lan\1n Z(s),P^{-1}Y(s)\1n\ran dW(s)\Big|\,\Big\}\\
\noalign{\smallskip}\displaystyle \,\les K\dbE\int_0^\infty\1n|\f(s)|^2ds
+2\dbE\Big\{\sup_{0\les t\les T}\Big|\int_t^T\2n\lan Z(s),P^{-1}Y(s)\ran dW(s)\Big|\,\Big\}\\
\noalign{\smallskip}\displaystyle \,\les{1\over4}\dbE\Big(\sup_{0\les t\les
T}\big|P^{-{1\over2}}Y(t)\big|^2\Big)
+K\dbE\int_0^\infty\1n|\f(s)|^2ds.\ea$$
Therefore (noting $P>0$ again),
\bel{E4}\dbE\Big(\sup_{0\les t\les T}|Y(t)|^2\Big)\les
K\dbE\int_0^\infty\1n|\f(s)|^2ds,\qq\forall T\in[0,\infty).\ee
Combining (\ref{E2}) and (\ref{E4}), making
use of Fatou's Lemma, yields (\ref{E}). \endpf

\ms

\bf Proposition 3.5. \sl Under the hypotheses of Proposition 3.4, we
have
\bel{Estimate}\dbE\int_0^\infty|Y(t)|^2dt\les K\dbE\int_0^\infty|\f(t)|^2dt.\ee

\ms

\it Proof. \rm Let $P>0$ be the matrix in the proof of Proposition
3.4. Applying It\^o's formula to $s\mapsto \lan
P^{-1}Y(s),Y(s)\ran$, one has that for any $0\les t<\infty$,
$$\ba{ll}
\noalign{\smallskip}\displaystyle  \dbE\lan P^{-1}Y(t),Y(t)\ran-\dbE\lan P^{-1}Y(0),Y(0)\ran\\
\noalign{\smallskip}\displaystyle \,=\dbE\int_0^t\1n\Big\{\1n-\1n\lan\1n P^{-1}\big[A^TY\1n+\1nC^TZ\1n+\1n\f\big],Y\1n\ran\1n
-\1n\lan\1n P^{-1}Y,A^TY\1n+\1nC^TZ\1n+\1n\f\1n\ran\1n+\1n\lan\1n P^{-1}Z,Z\1n\ran\Big\}ds\\
\noalign{\smallskip}\displaystyle \,=\dbE\int_0^t\1n\Big\{\1n-\1n\lan\1n PAP^{-1}Y,P^{-1}Y\1n\ran\1n-\1n\lan\1n A^TPP^{-1}Y,P^{-1}Y\1n\ran
\1n-2\lan \1nC^TZ\1n+\1n\f,P^{-1}Y\1n\ran\1n+\1n\lan\1n P^{-1}Z,Z\1n\ran\Big\}ds\\
\noalign{\smallskip}\displaystyle \,\ges\dbE\int_0^t\1n\Big\{\1n-\1n\lan\big[PA\1n+\1nA^TP\big]
P^{-1}Y,P^{-1}Y\1n\ran\1n-2\lan \1nC^TZ\1n+\1n\f,P^{-1}Y\1n\ran\Big\}ds.\ea$$
Let $\mu>0$ be the smallest eigenvalue of $-(PA+A^TP)>0$. By
Cauchy--Schwarz's inequality, we have
\bel{Estimate1}\ba{ll}
\noalign{\smallskip}\displaystyle \dbE\lan P^{-1}Y(t),Y(t)\ran-\dbE\lan P^{-1}Y(0),Y(0)\ran\\
\noalign{\smallskip}\displaystyle \,\ges\dbE\int_0^t\Big\{\mu|P^{-1}Y(s)|^2-{\mu\over2}|P^{-1}Y(s)|^2
-{4\over\mu}|C^TZ(s)|^2-{4\over\mu}|\f(s)|^2\Big\}ds,\q\forall t\in[0,\infty).\ea\ee
Letting $t\to\infty$ in (\ref{Estimate1}), one has
$$\dbE\lan P^{-1}Y(0),Y(0)\ran+{\mu\over2}\dbE\int_0^\infty\1n|P^{-1}Y(s)|^2ds
\les{4\over\mu}\dbE\int_0^\infty\1n\Big(|C^TZ(s)|^2+|\f(s)|^2\Big)ds.$$
Combining the a priori estimate (\ref{E}) we obtain the desired
estimate (\ref{Estimate}). \endpf

\ms

\it Proof of Theorem 3.3. \rm The uniqueness is an immediate
consequence of the a priori estimate (\ref{E}). We now prove the
existence. For $k=1,2,\cdots$, we set
$$\f_k(t)\deq1_{[0,k]}(t)\f(t),\q t\in[0,\infty).$$
Clearly, $\{\f_k(\cd)\}_{k=1}^\infty$ converges to $\f(\cd)$ in $L_\dbF^2(\dbR^n)$.

\ms

We now consider, for each $k$, the $L^2$-stable adapted solution
$(Y_k(\cd),Z_k(\cd))$ of the following BSDE:
\bel{BSDEk}dY_k(t)=-\big[A^TY_k(t)+C^TZ_k(t)+\f_k(t)\big]dt+Z_k(t)dW(t),\q t\in[0,\infty).\ee
The above can be solved as follows: on $[0,k]$,
$(Y_k(\cd),Z_k(\cd))$ is the adapted solution to the following
BSDE:
$$\left\{\2n\ba{ll}
\noalign{\smallskip}\displaystyle  dY_k(t)=-\big[A^TY_k(t)+C^TZ_k(t)+\f_k(t)\big]dt+Z_k(t)dW(t),\q t\in[0,k],\\
\noalign{\smallskip}\displaystyle  Y_k(k)=0,\ea\right.$$
and on $(k,\infty)$, it is identically equal to zero. By
Propositions 3.4 and 3.5, we have
$$\ba{ll}
\noalign{\smallskip}\displaystyle \dbE\Big(\sup_{0\les
t<\infty}\big|Y_k(t)-Y_j(t)\big|^2\Big)+\dbE\int_0^\infty\1n|Y_k(t)-Y_j(t)|^2dt
+\dbE\int_0^\infty\1n\big|Z_k(t)-Z_j(t)\big|^2dt\\
\noalign{\smallskip}\displaystyle \,\les K\dbE\int_0^\infty\1n\big|\f_k(t)-\f_j(t)\big|^2dt,\qq
\forall k, j.\ea$$
Therefore, there exists a $(Y(\cd),Z(\cd))\in\cX[0,\infty)\times L_\dbF^2(\dbR^n)$ such that
$$\dbE\Big(\sup_{0\les t<\infty}\big|Y_k(t)-Y(t)\big|^2\Big)+\dbE\int_0^\infty\big|Z_k(t)-Z(t)\big|^2dt
\to 0,\q\hb{as } k\to\infty,$$
which implies that $(Y(\cd),Z(\cd))$ is an $L^2$-stable adapted
solution of (\ref{BSDE}). \endpf

\section{Closed-Loop Optimal Controls}

In this section we discuss the closed-loop optimal controls of Problem (LQ). Let us first recall that for any $M\in\dbR^{m\times n}$, there exists a unique
matrix $M^\dag\in\dbR^{n\times m}$, called the (Moore-Penrose) {\it
pseudo-inverse} of $M$, satisfying the following (\cite{Penrose
1955}):
$$MM^\dag M=M,\q M^\dag MM^\dag=M^\dag,\q(MM^\dag)^T=MM^\dag,\q
(M^\dag M)^T=M^\dag M.$$
In addition, if $M\in\dbS^n$, then $M^\dag\in\dbS^n$, and
$$\ba{ll}
\noalign{\smallskip}\displaystyle  MM^\dag=M^\dag M;\qq M\ges0\iff M^\dag\ges0.\ea$$

\ms

\bf Lemma 4.1 (Extended Schur's Lemma \cite{Albert 1969}). \sl Let
$M\in\dbS^n$, $N\in\dbS^m$, $L\in\dbR^{n\times m}$. Then the
following conditions are equivalent:

\ms

{\rm(i)} $M-LN^\dag L^T\ges0$, $N\ges0$, and $L(I-NN^\dag)=0$.

\ms

{\rm(ii)} $\begin{pmatrix}M&L\\L^T&N\end{pmatrix}\ges0$.

\ms

\rm

Note that $L(I-NN^\dag)=0$ is equivalent to
$\cR(L^T)\subseteq\cR(N)$, where $\cR(\Lambda)$ is the range of a matrix
$\Lambda$. We now introduce the following notion.

\ms

\bf Definition 4.2. \rm A pair
$(\Th^*\1n,u^*(\cd))\1n\in\sS\times\1n L^2_\dbF(\dbR^m)$ is called a
{\it closed-loop optimal control} of Problem (LQ) if
\bel{closed-optimal}J(x;\Th^*X^*(\cd)+u^*(\cd))\les J(x;\Th
X(\cd)+u(\cd)),\q\forall(x,\Th,u(\cd))\in\dbR^n\1n\times\1n\sS\1n\times\1n
L^2_\dbF(\dbR^m).\ee

\ms

The following technical result, which is similar to Berkovitz's equivalence lemma for LQDG problems found in
\cite{Berkovitz 1971}, can be shown by a simple adaptation of \cite[Proposition 3.3]{Sun-Yong}.

\ms

\bf Proposition 4.3. \sl For
$(\Th^*\1n,u^*(\cd))\1n\in\1n\sS\1n\times\1n L^2_\dbF(\dbR^m)$, the
following statements are equivalent:

\ms

{\rm(i)} $(\Th^*,u^*(\cd))$ is a closed-loop optimal control of
Problem {\rm(LQ)}.

\ms

{\rm(ii)} For any $x\in\dbR^n$, and $u(\cd)\in L^2_\dbF(\dbR^m)$, the following holds:
\bel{}J(x;\Th^*X^*(\cd)+u^*(\cd))\les J(x;\Th^*X(\cd)+u(\cd)).\ee

\ms

\rm Now we present a characterization of closed-loop optimal controls of Problem (LQ) in terms of infinite horizon forward-backward
stochastic differential equations (FBSDE, for short).

\ms

\bf Theorem 4.4. \sl A pair $(\Th^*,u^*(\cd))\in\sS\1n\times\1n
L^2_\dbF(\dbR^m)$ is a closed-loop optimal control of Problem
{\rm(LQ)} if and only if for any $x\in\dbR^n$, the following FBSDE
admits an adapted solution
$(X^*(\cd),Y^*(\cd),Z^*(\cd))\in\cX[0,\infty)\times\cX[0,\infty)\times
L_\dbF^2(\dbR^n)$:
\bel{FBSDE}\left\{\2n\ba{ll}
\noalign{\smallskip}\displaystyle  dX^*(t)=\big\{(A+B\Th^*)X^*+Bu^*+b\big\}dt+\big\{(C+D\Th^*)X^*+Du^*+\si\big\}dW(t), \q t\ges0,\\
\noalign{\smallskip}\displaystyle  dY^*(t)=-\big\{A^TY^*+C^TZ^*+(Q+S^T\Th^*)X^*+S^Tu^*+q\big\}dt+Z^*dW(t),\q t\ges0,\\
\noalign{\smallskip}\displaystyle  X^*(0)=x,\ea\right.\ee
such that the following {\it stationarity condition} holds:
\bel{stationarity}Ru^*+B^TY^*+D^TZ^*+(S+R\Th^*)X^*+\rho=0,\q
\ae~\as\ee
and
\bel{convex}\dbE\int_0^\infty\1n\lan\begin{pmatrix}Q&S^T\\
                                                  S&R\end{pmatrix}
                                   \begin{pmatrix}X_0\\
                                                  \Th^*X_0+u\end{pmatrix},
                                   \begin{pmatrix}X_0\\
                                                  \Th^*X_0+u\end{pmatrix}\ran dt\ges0,
\q \forall u(\cd)\in L_\dbF^2(\dbR^m),\ee
where $X_0(\cd)$ is the solution of
\bel{convex-state}\left\{\2n\ba{ll}
\noalign{\smallskip}\displaystyle
dX_0(t)=\big\{\big[A+B\Th^*\big]X_0(t)+Bu(t)\big\}dt+\big\{\big[C+D\Th^*\big]X_0(t)+Du(t)\big\}dW(t),\q t\ges 0,\\
\noalign{\smallskip}\displaystyle  X_0(0)=0.\ea\right.\ee

\ms

\it Proof. \rm Consider the state equation
$$\left\{\2n\ba{ll}
\noalign{\smallskip}\displaystyle
dX(t)=\big\{\big[A+B\Th^*\big]X(t)+Bu(t)+b(t)\big\}dt\\
\noalign{\smallskip}\displaystyle \qq\qq\q+\big\{\big[C+D\Th^*\big]X(t)+Du(t)+\si(t)\big\}dW(t),\q t\ges 0,\\
\noalign{\smallskip}\displaystyle  X(0)=x,\ea\right.$$
with the cost functional
$$\ba{ll}
\noalign{\smallskip}\displaystyle  \wt J(x;u(\cd))\equiv J(x;\Th^*X(\cd)+u(\cd))\\
\noalign{\smallskip}\displaystyle \,=\dbE\int_0^\infty\1n\Big[\lan\begin{pmatrix}Q&S^T\\
                                                  S&R\end{pmatrix}
                                   \begin{pmatrix}X\\
                                                  \Th^*X\1n+\1nu\end{pmatrix},
                                   \begin{pmatrix}X\\
                                                  \Th^*X\1n+\1nu\end{pmatrix}\ran
                                                  \1n+2\lan\begin{pmatrix}q\\
                                                    \rho\end{pmatrix},
                                   \begin{pmatrix}X\\
                                                  \Th^*X\1n+\1nu\end{pmatrix}\ran\Big]dt\\
\noalign{\smallskip}\displaystyle \,=\dbE\int_0^\infty\1n\Big[\lan\begin{pmatrix}\wt Q&\wt S^T\\
                                                  \wt S&R\end{pmatrix}
                                   \begin{pmatrix}X\\u\end{pmatrix},
                                   \begin{pmatrix}X\\u\end{pmatrix}\ran\1n
                                   +2\lan\begin{pmatrix}\wt q\\ \rho\end{pmatrix},
                                   \begin{pmatrix}X\\u\end{pmatrix}\ran\Big]dt,\ea$$
where
$$\wt Q\1n=\1nQ\1n+\1n(\Th^*)^TS\1n+\1nS^T\Th^*\1n+\1n(\Th^*)^TR\Th^*,\q \wt S\1n=\1nS\1n+\1n R\Th^*,\q
\wt q\1n=\1nq\1n+\1n(\Th^*)^T\rho.$$
By Proposition 4.3, $(\Th^*,u^*(\cd))$ is a closed-loop optimal
control of Problem {\rm(LQ)} if and only if for any $x\in\dbR^n$,
$u^*(\cd)$ is an open-loop optimal control for the problem with the
above state equation and cost functional. For any $u(\cd)\in
L_\dbF^2(\dbR^m)$ and $\e\in\dbR$, let $X^\e(\cd)$ be the solution
of
$$\left\{\2n\ba{ll}
\noalign{\smallskip}\displaystyle  dX^\e(t)=\big\{\big[A+B\Th^*\big]X^\e(t)+B\big[u^*(t)+\e u(t)\big]+b(t)\big\}dt\\
\noalign{\smallskip}\displaystyle \qq\qq\q+\big\{\big[C+D\Th^*\big]X^\e(t)+D\big[u^*(t)+\e u(t)\big]+\si(t)\big\}dW(t),\qq t\ges0,\\
\noalign{\smallskip}\displaystyle  X^\e(0)=x.\ea\right.$$
Thus, $X_0(\cd)\equiv{X^\e(\cd)-X^*(\cd)\over\e}$ is independent of $\e$ and satisfies (\ref{convex-state}). Then
$$\ba{ll}
\noalign{\smallskip}\displaystyle  \wt J(x;u^*(\cd)+\e u(\cd))-\wt J(x;u^*(\cd))\\
\noalign{\smallskip}\displaystyle \,=\e\dbE\int_0^\infty\Big[\lan\begin{pmatrix}\wt Q&\wt S^T\\
                                                  \wt S&R\end{pmatrix}
                                   \begin{pmatrix}2X^*(t)+\e X_0(t)\\
                                                  2u^*(t)+\e u(t)\end{pmatrix},
                                   \begin{pmatrix}X_0(t)\\
                                                  u(t)\end{pmatrix}\ran
                                                  +2\lan\begin{pmatrix} \wt q(t)\\
                                                    \rho(t)\end{pmatrix},
                                   \begin{pmatrix}X_0(t)\\
                                                  u(t)\end{pmatrix}\ran\Big]dt\\
\noalign{\smallskip}\displaystyle \,=2\e\dbE\int_0^\infty\Big[\lan\wt QX^*,X_0\ran+\lan\wt SX^*,u\ran+\lan\wt SX_0,u^*\ran
+\lan Ru^*,u\ran+\lan\wt q,X_0\ran+\lan\rho,u\ran\Big]dt\\
\noalign{\smallskip}\displaystyle \q\,+\e^2\dbE\int_0^\infty\Big[\lan\wt QX_0(t),X_0(t)\ran
+2\lan\wt SX_0(t),u(t)\ran+\lan Ru(t),u(t)\ran\Big]dt\\
\noalign{\smallskip}\displaystyle \,=2\e\dbE\int_0^\infty\Big[\lan\wt QX^*+\wt S^Tu^*+\wt
q,X_0\ran+\lan\wt SX^*+Ru^*+\rho,u\ran\Big]dt\\
\noalign{\smallskip}\displaystyle \q\,+\e^2\dbE\int_0^\infty\Big[\lan\wt QX_0(t),X_0(t)\ran
+2\lan\wt SX_0(t),u(t)\ran+\lan Ru(t),u(t)\ran\Big]dt.\ea$$
Since $[A+B\Th^*,C+D\Th^*]$ is $L^2$-stable, by Theorem 3.3, the following BSDE:
$$\ba{ll}
\noalign{\smallskip}\displaystyle  dY^*\1n=\1n-\big\{(A\1n+\1nB\Th^*)^TY^*\1n+(C\1n+\1nD\Th^*)^TZ^*\1n+\1n\wt QX^*\1n
+\1n\wt S^Tu^*\1n+\1n\wt q\,\big\}dt\1n+\1nZ^*dW(t)\\
\noalign{\smallskip}\displaystyle \qq=\1n-\big\{A^TY^*\1n+\1nC^TZ^*\1n+\1nQX^*\1n+\1nS^T(\Th^*X^*\1n+\1nu^*)\1n+\1nq\\
\noalign{\smallskip}\displaystyle \qq\qq\,+(\Th^*)^T\big[B^TY^*\1n+\1nD^TZ^*\1n+\1n(S+R\Th^*)X^*\1n+\1n
Ru^*\1n+\1n\rho\big]\big\}dt\1n+\1nZ^*dW(t),\q t\ges 0\ea$$
admits a unique $L^2$-stable adapted solution $(Y^*(\cd),Z^*(\cd))$.
By It\^o's formula, we have
\bel{ttt}\ba{ll}
\noalign{\smallskip}\displaystyle \dbE\lan\1n
Y^*(t),X_0(t)\1n\ran\1n=\1n\dbE\int_0^t\1n\Big[\1n-\1n\lan(A\1n+\1nB\Th^*)^TY^*\1n+\1n(C\1n+\1nD\Th^*)^TZ^*\1n
+\1n\wt QX^*\1n+\1n\wt S^Tu^*\1n+\1n\wt q\,),X_0\1n\ran\\
\noalign{\smallskip}\displaystyle \qq\qq\qq\qq\qq\q+\1n\lan
\1nY^*,(A\1n+\1nB\Th^*)X_0\1n+\1nBu\1n\ran\1n+\lan\1n
Z^*,(C\1n+\1nD\Th^*)X_0\1n+\1nDu\1n\ran\1n\Big]ds\\
\noalign{\smallskip}\displaystyle \qq\qq\qq\q=\1n\dbE\int_0^t\1n\Big[\1n-\1n\lan\wt
QX^*\1n+\1n\wt S^Tu^*\1n+\1n\wt q,X_0\1n\ran+\1n\lan B^T\1n
Y^*\1n+\1n D^TZ^*,u\1n\ran\1n\Big]ds, \q \forall t\ges0.\ea\ee
Note that
$$\lim_{t\to\infty}|\,\dbE\lan Y^*(t),X_0(t)\ran|^2\les\lim_{t\to\infty}\dbE|Y^*(t)|^2\dbE|X_0(t)|^2=0.$$
Letting $t\to\infty$ in (\ref{ttt}), one has
$$\dbE\int_0^\infty\1n\lan\wt QX^*\1n+\1n\wt
S^Tu^*\1n+\1n\wt q,X_0\1n\ran ds=\dbE\int_0^\infty\1n\lan B^T
Y^*\1n+\1n D^TZ^*,u\1n\ran ds.$$
Hence,
$$\ba{ll}
\noalign{\smallskip}\displaystyle  \wt J(x;u^*(\cd)+\e u(\cd))-\wt J(x;u^*(\cd))\\
\noalign{\smallskip}\displaystyle \,=2\e\dbE\int_0^\infty\1n\Big[\lan\wt QX^*+\wt S^Tu^*+\wt
q,X_0\ran+\lan\wt SX^*+Ru^*+\rho,u\ran\Big]dt\\
\noalign{\smallskip}\displaystyle \q\,+\e^2\dbE\int_0^\infty\1n\Big[\lan\wt QX_0(t),X_0(t)\ran
+2\lan\wt SX_0(t),u(t)\ran+\lan Ru(t),u(t)\ran\Big]dt\\
\noalign{\smallskip}\displaystyle \,=2\e\dbE\int_0^\infty\1n\lan B^TY^*+D^TZ^*+\wt SX^*+Ru^*+\rho,u\ran dt\\
\noalign{\smallskip}\displaystyle \q\,+\e^2\dbE\int_0^\infty\1n\lan\begin{pmatrix}Q&S^T\\
                                                  S&R\end{pmatrix}
                                   \begin{pmatrix}X_0\\
                                                  \Th^*X_0+u\end{pmatrix},
                                   \begin{pmatrix}X_0\\
                                                  \Th^*X_0+u\end{pmatrix}\ran dt.\ea$$
Therefore, $(\Th^*,u^*(\cd))$ is a closed-loop optimal control of
Problem {\rm(LQ)} if and only if (\ref{stationarity}) and
(\ref{convex}) hold. Consequently, $(Y^*(\cd),Z^*(\cd))$ solves the
following BSDE:
$$dY^*=-\big\{A^TY^*+C^TZ^*+QX^*+S^T(\Th^*X^*+u^*)+q\big\}dt+Z^*dW(t),\q t\ges0.$$
This completes the proof. \endpf

\ms

As a consequence, we have the following result.

\ms

\bf Corollary 4.5. \sl If $(\Th^*,u^*(\cd))$ is a closed-loop
optimal control of Problem {\rm(LQ)}, then $(\Th^*,0)$ is a
closed-loop optimal control of Problem $\hb{\rm{(LQ)}}^0$.

\ms

\it Proof. \rm Let $(\Th^*,u^*(\cd))$ be a closed-loop optimal
control of Problem (LQ). Then, by Theorem 4.4, (\ref{convex})
holds, and for any $x\in\dbR^n$, FBSDE (\ref{FBSDE}) admits an
adapted solution
$(X^*(\cd),Y^*(\cd),Z^*(\cd))\1n\in\1n\cX[0,\infty)\times\1n\cX[0,\infty)\1n\times\1n
L_\dbF^2(\dbR^n)$ satisfying (\ref{stationarity}). Since FBSDE
(\ref{FBSDE}) admits a solution for each $x\in\dbR^n$, and
$(\Th^*,u^*(\cd))$ is independent of $x$, by subtracting solutions
corresponding $x$ and $0$, the later from the former, we see that
for any $x\in\dbR^n$, the following FBSDE:
$$\left\{\2n\ba{ll}
\noalign{\smallskip}\displaystyle  dX=(A+B\Th^*)Xdt+(C+D\Th^*)XdW(t),\q t\ges0,\\
\noalign{\smallskip}\displaystyle  dY=-\big[A^TY+C^TZ+(Q+S^T\Th^*)X\big]dt+ZdW(t),\q t\ges0,\\
\noalign{\smallskip}\displaystyle  X(0)=x,\ea\right.$$
admits an adapted solution
$(X(\cd),Y(\cd),Z(\cd))\1n\in\1n\cX[0,\infty)\1n\times\1n\cX[0,\infty)\1n\times\1n L_\dbF^2(\dbR^n)$ satisfying
$$B^TY+D^TZ+(S+R\Th^*)X=0,\q \ae~\as$$
Again, by Theorem 4.4, we see that $(\Th^*,0)$ is a closed-loop
optimal control of Problem $\hb{(LQ)}^0$. \endpf

\ms

The following theorem gives a necessary condition for the existence of a closed-loop optimal control of Problem {\rm(LQ)}.

\ms

\bf Theorem 4.6. \sl Suppose Problem {\rm(LQ)} admits a closed-loop optimal control. Then the following ARE:
\bel{Riccati}
PA\1n+\1nA^TP\1n+\1nC^TPC\1n+\1nQ\1n
-\1n\big(PB\1n+\1nC^TPD\1n+\1nS^T\big)(R\1n+\1nD^TPD)^\dag\big(B^TP\1n+\1nD^TPC\1n+\1nS\big)\1n=\1n0\ee
admits a solution $P\in\dbS^n$ such that
\bel{R-C}R+D^TPD\ges0,\qq \cR\big(B^TP+D^TPC+S\big)\subseteq\cR\big(R+D^TPD\big),\ee
and there exists a $\Pi\in\dbR^{m\times n}$ such that
\bel{Th}-(R+D^TPD)^\dag(B^TP+D^TPC+S)+\big[I-(R+D^TPD)^\dag(R+D^TPD)\big]\Pi\ee
is a stabilizer of $[A,C;B,D]$.

\ms

\it Proof. \rm Let $(\Th^*,u^*(\cd))$ be a closed-loop optimal
control of Problem (LQ). Then, by Corollary 4.5, $(\Th^*,0)$ is a
closed-loop optimal control of Problem $\hb{(LQ)}^0$, and hence
Problem $\hb{(LQ)}^0$ is finite. Lemma 2.7 yields that the set $\sP$
has a maximal element $P\in\sP$ such that $V^0(x)=\lan Px,x\ran$, and
\bel{matrix form}\begin{pmatrix}\sM(P)&\sL(P)\\\sL(P)^T&\sN(P)\end{pmatrix}\ges0.\ee
Applying Lemma 4.1 to (\ref{matrix form}), we have
\bel{N-1}\sM(P)-\sL(P)\sN(P)^\dag\sL(P)^T\ges0,\ee
\bel{N-2}\sN(P)\ges0,\qq \sL(P)\big[I-\sN(P)\sN(P)^\dag\big]=0.\ee
Note that (\ref{N-2}) is equivalent to (\ref{R-C}). Let $X^*(\cd)$ be the solution of
$$\left\{\2n\ba{ll}
\noalign{\smallskip}\displaystyle  dX^*(t)=\big[A+B\Th^*\big]X^*(t)dt+\big[C+D\Th^*\big]X^*(t)dW(t),\q t\ges 0,\\
\noalign{\smallskip}\displaystyle  X(0)=x.\ea\right.$$
Applying It\^o's formula to $t\to \lan PX(t),X(t)\ran$, one has
$$\ba{ll}
\noalign{\smallskip}\displaystyle  \lan Px,x\ran\1n=\1n-\dbE\int_0^\infty\2n\Big\{\1n\lan\big[P(A\1n+\1n
B\Th^*)\1n+\1n(A\1n+\1nB\Th^*)^TP\big]X,X\ran\1n
+\1n\lan P(C\1n+\1nD\Th^*)X,(C\1n+\1nD\Th^*)X\ran\1n\Big\}dt\\
\noalign{\smallskip}\displaystyle \qq\qq\2n=\1n-\dbE\int_0^\infty\2n\lan\big[(PA\1n+\1nA^TP\1n+\1nC^TPC)\1n+\1n(PB\1n+\1nC^TPD)\Th^*\\
\noalign{\smallskip}\displaystyle \qq\qq\qq\qq\q+(\Th^*)^T(B^TP\1n+\1nD^TPC)\1n+\1n(\Th^*)^TD^TPD\Th^*\big]X,X\ran dt\\
\noalign{\smallskip}\displaystyle \qq\qq\2n=\1n-\dbE\int_0^\infty\2n\lan\big[\sM(P)\1n+\1n\sL(P)\Th^*\1n+\1n(\Th^*)^T\sL(P)^T\1n
+\1n(\Th^*)^T\sN(P)\Th^*\big]X,X\ran dt\\
\noalign{\smallskip}\displaystyle \qq\qq\q\2n+\dbE\int_0^\infty\2n\lan\big[Q\1n+\1nS^T\Th^*\1n+\1n(\Th^*)^TS\1n+\1n(\Th^*)^TR\Th^*\big]X,X\ran
dt.\ea$$
Then we have (noting (\ref{N-2}))
\bel{N-3}\ba{ll}
\noalign{\smallskip}\displaystyle
V^0(x)\1n=\1nJ^0(x,\Th^*X(\cd))=\dbE\int_0^\infty\2n\lan\big[Q\1n+\1nS^T\Th^*\1n+\1n(\Th^*)^TS\1n
+\1n(\Th^*)^TR\Th^*\big]X,X\ran dt\\
\noalign{\smallskip}\displaystyle \qq\q\2n=\1n\lan Px,x\ran\1n+\dbE\int_0^\infty\2n\lan
\big[\sM(P)\1n+\1n\sL(P)\Th^*\1n+\1n(\Th^*)^T\sL(P)^T\1n
+\1n(\Th^*)^T\sN(P)\Th^*\big]X,X\ran dt\\
\noalign{\smallskip}\displaystyle \qq\q\2n=\1n\lan Px,x\ran\1n+\dbE\int_0^\infty\2n\lan \big[\sM(P)\1n-\1n\sL(P)\sN(P)^\dag\sL(P)^T\,\big]X,
X\ran dt\\
\noalign{\smallskip}\displaystyle \qq\qq\qq\q+\dbE\int_0^\infty\2n\lan\sN(P)\big[\Th^*\1n+\1n\sN(P)^\dag\sL(P)^T\,\big]X,
[\Th^*\1n+\1n\sN(P)^\dag\sL(P)^T\,\big]X
\ran dt.\ea\ee
Due to the equality $V^0(x)=\lan Px,x\ran$ and
(\ref{N-1})--(\ref{N-3}), each of the two integrands on the
right-hand side of (\ref{N-3}) must be zero almost everywhere. Hence, we obtain
$$\sM(P)-\sL(P)\sN(P)^\dag\sL(P)^T=0,$$
that is, $P$ is a solution of (\ref{Riccati}), and
$$\sN(P)^{1\over2}\big[\Th^*+\sN(P)^\dag\sL(P)^T\big]=0,$$
which, together with (\ref{N-2}), gives
\bel{N-4}\sN(P)\Th^*+\sL(P)^T=0.\ee
Since $\sN(P)\sN(P)^\dag$ is an orthogonal projection, we have
$$\Th^*=-\sN(P)^\dag\sL(P)^T+\big[I-\sN(P)^\dag\sN(P)\big]\Pi\in\sS,$$
for some $\Pi\in\dbR^{n\times m}$. \endpf

\ms

We point out that the sufficiency of the above result can also be
stated and proved, which is a special case of the corresponding
result for two-person zero-sum differential games (see the next
section). Hence, to avoid a repeating presentation, we prefer not to
give the details here.

\ms

\section{Open-Loop and Closed-Loop Saddle Points}

We now return to our differential games. For notational simplicity, we let $m=m_1+m_2$ and denote
$$\ba{ll}
\noalign{\smallskip}\displaystyle  B=(B_1,B_2),\q D=(D_1,D_2),\\
\noalign{\smallskip}\displaystyle  S=\begin{pmatrix}S_1\\ S_2\end{pmatrix},\q R=\begin{pmatrix}R_{11}&R_{12}\\
R_{21}&R_{22}\end{pmatrix}\equiv\begin{pmatrix}R_1\\ R_2\end{pmatrix},\q\rho(\cd)=\begin{pmatrix}\rho_1(\cd)\\
\rho_2(\cd)\end{pmatrix},\q u(\cd)=\begin{pmatrix}u_1(\cd)\\u_2(\cd)\end{pmatrix}.\ea$$
With such notations, the state equation becomes
\bel{state5.1}\left\{\2n\ba{ll}
\noalign{\smallskip}\displaystyle
dX(t)\1n=\1n\big[A(t)X(t)\1n+\1nB(t)u(t)\1n+\1nb(t)\big]dt\1n+\1n\big[C(t)X(t)\1n+\1nD(t)u(t)\1n+\1n\si(t)\big]dW(t),\q t\ges0,\\
\noalign{\smallskip}\displaystyle  X(0)\1n=\1nx,\ea\right.\ee
and the performance functional becomes
\bel{cost5.2}
J(x;u_1(\cd),u_2(\cd))\1n=\1nJ(x;u(\cd))\1n=\1n\dbE\1n\int_0^\infty\1n\Big[\1n\lan\1n\begin{pmatrix}Q&S^T\\
                                          S&R\end{pmatrix}\begin{pmatrix}
                                          X(t)\\ u(t)\end{pmatrix},
                                          \begin{pmatrix}
                                          X(t)\\u(t)\end{pmatrix}\1n\ran\1n+
                                          2\lan\1n\begin{pmatrix}q(t)\\ \rho(t)
                                          \end{pmatrix},\begin{pmatrix}X(t)\\
                                          u(t)\
                                          \end{pmatrix}\1n\ran\1n\Big] dt.\ee
Also, when $b(\cd),\si(\cd),q(\cd),\rho(\cd)\1n=\1n0$, we denote the corresponding Problem (LQG)
by Problem $\hb{(LQG)}^0$ and the corresponding performance functional by
$J^0(x;u_1(\cd),u_2(\cd))$. Similar to Problem (LQ), we will assume
(H1) for the system $[A,C;B,D]$, and we also denote
$$\sM(P)\1n=\1nPA\1n+\1nA^TP\1n+\1nC^TPC\1n+\1nQ,\q \sL(P)\1n=\1nPB\1n+\1nC^TPD\1n+\1nS^T,
\q \sN(P)\1n=\1nR\1n+\1nD^TPD;\qq \forall P\in\dbS^n.$$
Moreover, for
$\Th_i\in\dbR^{m_i\times n}$, $i=1,2$, we let
$$\ba{ll}
\noalign{\smallskip}\displaystyle \sS_1(\Th_2)=\Big\{\Th_1\in\dbR^{m_1\times
n}\bigm|(\Th_1^T,\Th_2^T)^T \hb{ is a stabilizer of }
[A,C;B,D]\Big\},\\
\noalign{\smallskip}\displaystyle \sS_2(\Th_1)=\Big\{\Th_2\in\dbR^{m_2\times
n}\bigm|(\Th_1^T,\Th_2^T)^T \hb{ is a stabilizer of }
[A,C;B,D]\Big\}.\ea$$
Note that in general, say, $\sS_1(\Th_2)$ is not necessarily
non-empty for some $\Th_2\in\dbR^{m_2\times n}$. However, if
$\Th\equiv(\Th_1^T,\Th_2^T)^T\in\sS[A,C;B,D]$, then both
$\sS_1(\Th_2)$ and $\sS_2(\Th_1)$ are non-empty. Also, for any
$x\in\dbR^n$, we let $\cU_{ad}(x)$ be the set of all
$u(\cd)\equiv(u_1(\cd),u_2(\cd))\in L^2_\dbF(\dbR^m)$ such that the
corresponding state $X(\cd)\equiv
X(\cd\,;x,u(\cd))\in\cX[0,\infty)$.

\ms

\bf Definition 5.1. \rm For any given $x\in\dbR^n$, a pair $(\bar
u_1(\cd),\bar u_2(\cd))\in\cU_{ad}(x)$ is called an {\it open-loop
saddle point} of Problem (LQG) if
\bel{saddle-open}\ba{ll}
\noalign{\smallskip}\displaystyle  J(x;\bar u_1(\cd),u_2(\cd))\les J(x;\bar u_1(\cd),\bar
u_2(\cd))\les J(x;u_1(\cd),\bar u_2(\cd)),\ea\ee
for any $(u_1(\cd),u_2(\cd))\in L^2_\dbF(\dbR^m)$ such that
$J(x;\bar u_1(\cd),u_2(\cd))$ and $J(x;u_1(\cd),\bar u_2(\cd))$ are
defined.

\ms

\bf Definition 5.2. \rm A 4-tuple
$(\Th_1^*,u_1^*(\cd);\Th_2^*,u_2^*(\cd))\1n\in\dbR^{m_1\1n\times
n}\1n\times\1n L^2_\dbF(\dbR^{m_1})\1n\times\1n\dbR^{m_2\1n\times
n}\1n\times\1n L^2_\dbF(\dbR^{m_2})$ is called a {\it closed-loop
saddle point} of Problem (LQG) if

\ms

{\rm(i)} $\Th^*\equiv((\Th_1^*)^T,(\Th_2^*)^T)^T\in\sS[A,C;B,D]$,

\ms

{\rm(ii)} for any $x\in\dbR^n$,
$(\Th_1,\Th_2)\1n\in\sS_1(\Th_2^*)\1n\times\1n\sS_2(\Th_1^*)$ and
$(u_1(\cd),u_2(\cd))\1n\in L^2_\dbF(\dbR^{m_1})\1n\times\1n
L^2_\dbF(\dbR^{m_2})$,
\bel{saddle-closed}\ba{ll}
\noalign{\smallskip}\displaystyle  J(x;\Th_1^*X(\cd)+u_1^*(\cd),\Th_2X(\cd)+u_2(\cd))\les
J(x;\Th_1^*X^*(\cd)+u_1^*(\cd),
\Th_2^*X^*(\cd)+u_2^*(\cd))\\
\noalign{\smallskip}\displaystyle \qq\qq\qq\qq\qq\qq\qq\qq\,\les
J(x;\Th_1X(\cd)+u_1(\cd),\Th_2^*(\cd)X(\cd)+u_2^*(\cd)).\ea\ee

\ms

\bf Remark 5.3. \rm (a) Although both players are non-cooperative,
when choosing $\Th_i$ $(i=1,2)$, they prefer to at least work
together so that $\Th=((\Th_1)^T,(\Th_2)^T)^T$ is a stabilizer of
$[A,C;B,D]$ (and the system will not be crashed). Thus, in
Definition 5.2, we only require $\Th^*$ being a stabilizer of
$[A,C;B,D]$ rather than $\Th_i^*$ being a a stabilizer of
$[A,C;B_i,D_i]$.

\ms

(b) By a similar method used in \cite{Sun-Yong}, one can show that
condition (ii) in Definition 5.2 is equivalent to the following:

\ms

(ii)$'$ for any $x\in\dbR^n$ and $(u_1(\cd),u_2(\cd))\1n\in
L^2_\dbF(\dbR^{m_1})\1n\times\1n L^2_\dbF(\dbR^{m_2})$,
\bel{saddle-closed}\ba{ll}
\noalign{\smallskip}\displaystyle  J(x;\Th_1^*X(\cd)+u_1^*(\cd),\Th_2^*X(\cd)+u_2(\cd))\les J(x;\Th_1^*X^*(\cd)+u_1^*(\cd),
\Th_2^*X^*(\cd)+u_2^*(\cd))\\
\noalign{\smallskip}\displaystyle \qq\qq\qq\qq\qq\qq\qq\qq\,\les
J(x;\Th_1^*X(\cd)+u_1(\cd),\Th_2^*(\cd)X(\cd)+u_2^*(\cd)).\ea\ee

\ms

Let $\Th^*\1n=\1n((\Th_1^*)^T,(\Th_2^*)^T)^T\in\sS[A,C;B,D]$ and
$u^*(\cd)=(u_1^*(\cd)^T,u_2^*(\cd)^T)^T\1n\in L^2_\dbF(\dbR^m)$. We
look at the following state equation:
$$\left\{\2n\ba{ll}
\noalign{\smallskip}\displaystyle
dX(t)\1n=\1n\Big\{\1n\big[A\1n+\1nB\Th^*\big]X(t)\1n+\1nBu(t)\1n+\1nb(t)\1n\Big\}dt\1n
+\1n\Big\{\1n\big[C\1n+\1nD\Th^*\big]X(t)\1n+\1nDu(t)\1n+\1n\si(t)\Big\}dW(t),\q t\ges 0,\\
\noalign{\smallskip}\displaystyle  X(0)\1n=\1nx,\ea\right.$$
and the following performance functional:
$$\ba{ll}
\noalign{\smallskip}\displaystyle \wt J(x;u_1(\cd),u_2(\cd))\equiv J(x;\Th_1^*X(\cd)+u_1(\cd),\Th_2^*X(\cd)+u_2(\cd))\\
\noalign{\smallskip}\displaystyle \qq\qq\qq\q\ =\dbE\int_0^\infty\1n\Big[\lan\begin{pmatrix}\wt Q&\wt S^T\\
                                                  \wt S&R\end{pmatrix}
                                   \begin{pmatrix}X\\u\end{pmatrix},
                                   \begin{pmatrix}X\\u\end{pmatrix}\ran\1n
                                   +2\lan\begin{pmatrix}\wt q\\ \rho\end{pmatrix},
                                   \begin{pmatrix}X\\u\end{pmatrix}\ran\Big]dt,\ea$$
where
$$\wt Q\1n=\1nQ\1n+\1n(\Th^*)^TS\1n+\1nS^T\Th^*\1n+\1n(\Th^*)^TR\Th^*,\q \wt S\1n=\1nS\1n+\1nR\Th^*,\q
\wt q\1n=\1nq\1n+\1n(\Th^*)^T\rho.$$
From (ii)$'$ of Remark 5.3, we see that
$(\Th_1^*,u_1^*(\cd);\Th_2^*,u_2^*(\cd))$ is a closed-loop saddle
point of Problem {\rm(LQG)} if and only if $(u_1^*(\cd),u_2^*(\cd))$
is an open-loop saddle point  for the problem with the above state
equation and performance functional. Applying the ideal used in the
proof of Theorem 4.4 (see also \cite[Theorem 4.1]{Sun-Yong}), we see
that $(\Th_1^*,u_1^*(\cd);\Th_2^*,u_2^*(\cd))$ is a closed-loop
saddle point of Problem {\rm(LQG)} if and only if for any
$x\in\dbR^n$, the adapted solution
$(X^*(\cd),Y^*(\cd),Z^*(\cd))\1n\in\1n\cX[0,\infty)\1n\times\1n\cX[0,\infty)\1n\times\1n
L_\dbF^2(\dbR^n)$ of the following FBSDE:
\bel{FBSDE-G}\left\{\2n\ba{ll}
\noalign{\smallskip}\displaystyle  dX^*(t)\1n=\1n\big\{(A\1n+\1nB\Th^*)X^*\1n+\1nBu^*\1n+\1nb\big\}dt+\big\{(C\1n+\1nD\Th^*)X^*\1n+\1nDu^*\1n+\1n\si\big\}dW(t), \q t\ges0,\\
\noalign{\smallskip}\displaystyle  dY^*(t)\1n=\1n-\big\{(A\1n+\1nB\Th^*)^TY^*\1n+\1n(C\1n+\1nD\Th^*)^TZ^*\1n+\1n\wt QX^*\1n+\1n\wt S^Tu^*\1n+\1n\wt q\,\big\}dt+Z^*dW(t),\q t\ges0,\\
\noalign{\smallskip}\displaystyle  X^*(0)\1n=\1nx,\ea\right.\ee
satisfies the following stationarity condition:
\bel{stationarity-G}Ru^*+B^TY^*+D^TZ^*+\wt SX^*+\rho=0,\q
\ae~\as\ee
and the following convexity-concavity conditions hold: For $i=1,2$,
\bel{convex-concave}(-1)^{i-1}\dbE\1n\int_0^\infty\2n\lan\begin{pmatrix}\wt Q&\wt S_i^T\\
                                                  \wt S_i&R_{ii}\end{pmatrix}
                                   \begin{pmatrix}X_i\\
                                                  u_i\end{pmatrix},
                                   \begin{pmatrix}X_i\\
                                                  u_i\end{pmatrix}\ran dt\ges0,
\q \forall u_i(\cd)\in L_\dbF^2(\dbR^{m_i}),\ee
where $\wt S_i=S_i+R_i\Th^*$ and $X_i(\cd)$ is the solution of
\bel{Xi}\left\{\2n\ba{ll}
\noalign{\smallskip}\displaystyle
dX_i(t)=\Big\{\big[A+B\Th^*\big]X_i(t)+B_iu_i(t)\Big\}dt+\Big\{\big[C+D\Th^*\big]X_i(t)+D_iu_i(t)\Big\}dW(t),
\q t\ges 0,\\
\noalign{\smallskip}\displaystyle  X_i(0)=0.\ea\right.\ee
Applying the method used in the proof of Corollary 4.5, we obtain the following result.

\ms

\bf Proposition 5.4. \sl If
$(\Th_1^*,u_1^*(\cd);\Th_2^*,u_2^*(\cd))$ is a closed-loop saddle
point of Problem {\rm(LQG)}, then $(\Th_1^*,0;\Th_2^*,0)$ is a
closed-loop saddle point of Problem $\hb{\rm{(LQG)}}^0$.

\ms

\rm Next, we consider the following algebraic Riccati equation:
\bel{R-G}\left\{\2n\ba{ll}
\noalign{\smallskip}\displaystyle  PA\1n+\1nA^TP\1n+\1nC^TPC\1n+\1nQ\1n-
\1n\big(PB\1n+\1nC^TPD\1n+\1nS^T\big)(R\1n+\1nD^TPD)^\dag\big(B^TP\1n+\1nD^TPC\1n+\1nS\big)\1n=\1n0,\\
\noalign{\smallskip}\displaystyle \cR\big(B^TP+D^TPC+S\big)\subseteq\cR\big(R+D^TPD\big),\\
\noalign{\smallskip}\displaystyle  R_{11}+D_1^TPD_1\ges0,\q R_{22}+D_2^TPD_2\les0.\ea\right.\ee

\ms

\bf Definition 5.5. \rm A $P\in\dbS^n$ is called a {\it stabilizing
solution} of (\ref{R-G}) if $P$ is a solution to (\ref{R-G}) and
there exists a $\Pi\1n\in\1n\dbR^{m\times n}$ such that
$$-\sN(P)^\dag\sL(P)^T\1n+\1n\big[I\1n-\1n\sN(P)^\dag\sN(P)\big]\Pi
\in\sS[A,C;B,D].$$

\ms

Now we give a necessary condition for the existence of closed-loop saddle points of Problem $\hb{\rm{(LQG)}}^0$.

\ms

\bf Proposition 5.6. \sl Suppose Problem $\hb{\rm{(LQG)}}^0$ admits
a closed-loop saddle point. Then ARE (\ref{R-G}) admits a
stabilizing solution $P$.

\ms

\it Proof. \rm We assume without loss of generality that
$(\Th_1^*,0;\Th_2^*,0)$ is a closed-loop saddle point of Problem
$\hb{\rm{(LQG)}}^0$. Set
$$V^0(x)\deq J^0(x;\Th_1^*X^*(\cd),\Th_2^*X^*(\cd)).$$
It is easily seen that $V^0(\cd)$ is a quadratic form, that is,
there is a $P\in\dbS^n$ such that
$$V^0(x)=\lan Px,x\ran,\q \forall x\in\dbR^n.$$
Consider the state equation
$$\left\{\2n\ba{ll}
\noalign{\smallskip}\displaystyle
dX_1(t)\1n=\1n\Big\{\big[A\1n+\1nB_2\Th_2^*\big]X_1(t)\1n+\1nB_1u_1(t)\Big\}dt
\1n+\1n\Big\{\big[C\1n+\1nD_2\Th_2^*\big]X_1(t)\1n+\1nD_1u_1(t)\Big\}dW(t),\q t\ges 0,\\
\noalign{\smallskip}\displaystyle  X_1(0)\1n=\1nx,\ea\right.$$
with the cost functional
$$\ba{ll}
\noalign{\smallskip}\displaystyle  J_1(x;u_1(\cd))\equiv J^0(x;u_1(\cd),\Th_2^*X_1(\cd))
=\dbE\int_0^\infty\2n\lan\begin{pmatrix}Q&S_1^T&S_2^T\\
                                          S_1&R_{11}&R_{12}\\
                                          S_2&R_{21}&R_{22}\end{pmatrix}
                                          \begin{pmatrix}X_1\\ u_1\\ \Th_2^*X_1
                                          \end{pmatrix},\begin{pmatrix}X_1\\ u_1\\ \Th_2^*X_1
                                          \end{pmatrix}\ran dt\\
\noalign{\smallskip}\displaystyle \,=\dbE\int_0^\infty\2n\Big\{\1n\lan[Q\1n+\1n(\Th_2^*)^TR_{22}\Th_2^*\1n+\1n(\Th_2^*)^TS_2\1n+\1nS_2^T\Th_2^*]
X_1,X_1\ran\1n+\1n\lan
R_{11}u_1,u_1\ran\1n+2\lan(S_1\1n+\1nR_{12}\Th_2^*)X_1,u_1\ran\1n\Big\}dt.\ea$$
Then $(\Th_1^*,0)$ is a closed-loop optimal control of Problem
$\hb{\rm{(LQ)}}^0$ with the above state equation and cost
functional, and the value function of the above problem is given by $\lan Px,x\ran$. By Theorem 4.6, $P$ solves the following ARE:
\bel{ARE1}
P\wt A_1\1n+\1n\wt A_1^TP\1n+\1n\wt C_1^TP\wt C_1\1n+\1n\wt Q_1\1n-\1n\big(PB_1\1n+\1n\wt C_1^TPD_1\1n+\1n\wt S_1^T\big)(R_{11}\1n+\1nD_1^TPD_1)^\dag\big(B_1^TP\1n+\1nD_1^TP\wt C_1\1n+\1n\wt S_1\big)\1n=\1n0\ee
and (noting (\ref{N-4}))
\bel{ARE1.1}R_{11}\1n+\1nD_1^TPD_1\ges0,\q (R_{11}\1n+\1nD_1^TPD_1)\Th_1^*+\big(B_1^TP\1n+\1nD_1^TP\wt C_1\1n+\1n\wt S_1\big)=0,\ee
where
$$\wt A_1\1n=\1nA\1n+\1nB_2\Th_2^*,\q \wt C_1\1n=\1nC\1n+\1nD_2\Th_2^*,\q
\wt
Q_1\1n=\1nQ\1n+\1n(\Th_2^*)^TR_{22}\Th_2^*\1n+\1n(\Th_2^*)^TS_2\1n+\1nS_2^T\Th_2^*,\q
\wt S_1\1n=\1nS_1\1n+\1nR_{12}\Th_2^*.$$
Similarly, by considering the state equation
$$\left\{\2n\ba{ll}
\noalign{\smallskip}\displaystyle
dX_2(t)\1n=\1n\Big\{\big[A\1n+\1nB_1\Th_1^*\big]X_2(t)\1n+\1nB_2u_2(t)\Big\}dt\1n+\1n
\Big\{\big[C\1n+\1nD_1\Th_1^*\big]X_2(t)\1n+\1nD_2u_2(t)\Big\}dW(t),\q t\ges 0,\\
\noalign{\smallskip}\displaystyle  X_2(0)\1n=\1nx,\ea\right.$$
with the cost functional $J_2(x;u_2(\cd))\equiv
-J^0(x;\Th_1^*X_2(\cd),u_2(\cd))$, we have
\bel{ARE2.1}R_{22}\1n+\1nD_2^TPD_2\les0,\q
(R_{22}\1n+\1nD_2^TPD_2)\Th_2^*+\big(B_2^TP\1n+\1nD_2^TP\wt
C_2\1n+\1n\wt S_2\big)=0,\ee
where
$$\wt C_2\1n=\1nC\1n+\1nD_1\Th_1^*,\q \wt S_2\1n=\1nS_2\1n+\1nR_{21}\Th_1^*.$$
Let $\Th^*=((\Th^*_1)^T,(\Th^*_2)^T)^T$. Combining
(\ref{ARE1.1}) and (\ref{ARE2.1}), one has
\bel{ARE3}(R+D^TPD)\Th^*+\big(B^TP+D^TPC+S\big)=0,\ee
which implies
$$\cR\big(B^TP\1n+\1nD^TPC\1n+\1nS\big)\subseteq\cR\big(R\1n+\1nD^TPD\big).$$
Since $\sN(P)^\dag\sN(P)$ is an orthogonal projection, there exists a
$\Pi\in\dbR^{m\times n}$ such that
\bel{K^*}\Th^*=-\sN(P)^\dag\sL(P)^T+\big[I-\sN(P)^\dag\sN(P)\big]\Pi\in\sS[A,C;B,D].\ee
Using (\ref{ARE1})--(\ref{ARE3}), we have
\be{}\ba{ll}
\noalign{\smallskip}\displaystyle  0=P\wt A_1\1n+\1n\wt A_1^TP\1n+\1n\wt C_1^TP\wt C_1\1n+\1n\wt Q_1\1n-\1n\big(PB_1\1n+\1n\wt C_1^TPD_1\1n+\1n\wt S_1^T\big)(R_{11}\1n+\1nD_1^TPD_1)^\dag\big(B_1^TP\1n+\1nD_1^TP\wt C_1\1n+\1n\wt S_1\big)\\
\noalign{\smallskip}\displaystyle \q\1n=P\wt A_1\1n+\1n\wt A_1^TP\1n+\1n\wt C_1^TP\wt C_1\1n+\1n\wt Q_1\1n-\1n(\Th_1^*)^T(R_{11}\1n
+\1nD_1^TPD_1)\Th_1^*\\
\noalign{\smallskip}\displaystyle \q\1n=PA\1n+\1nA^TP\1n+\1nC^TPC\1n+\1nQ\1n+\1n(\Th_2^*)^T(R_{22}\1n+\1nD_2^TPD_2)\Th_2^*\1n
-\1n(\Th_1^*)^T(R_{11}\1n+\1nD_1^TPD_1)\Th_1^*\\
\noalign{\smallskip}\displaystyle \qq+\big(PB_2\1n+\1nC^TPD_2\1n+\1nS_2^T\big)\Th_2^*\1n+\1n(\Th_2^*)^T
\big(B_2^TP\1n+\1nD_2^TPC\1n+\1nS_2\big)\\
\noalign{\smallskip}\displaystyle \q\1n=PA\1n+\1nA^TP\1n+\1nC^TPC\1n+\1nQ\1n-\1n(\Th_1^*)^T(R_{11}\1n
+\1nD_1^TPD_1)\Th_1^*-\1n(\Th_2^*)^T(R_{22}\1n+\1nD_2^TPD_2)\Th_2^*\1n\\
\noalign{\smallskip}\displaystyle \qq+\big[(\Th_2^*)^T(R_{22}\1n+\1nD_2^TPD_2)\1n+\1n\big(PB_2\1n+\1nC^TPD_2\1n+\1nS_2^T\big)\big]\Th_2^*\\
\noalign{\smallskip}\displaystyle \qq+(\Th_2^*)^T\big[\big(B_2^TP\1n+\1nD_2^TPC\1n+\1nS_2\big)\1n+\1n(R_{22}\1n+\1nD_2^TPD_2)\Th_2^*\big]\\
\noalign{\smallskip}\displaystyle \q\1n=PA\1n+\1nA^TP\1n+\1nC^TPC\1n+\1nQ\1n-\1n(\Th_1^*)^T(R_{11}\1n+\1nD_1^TPD_1)\Th_1^*-\1n(\Th_2^*)^T(R_{22}\1n+\1nD_2^TPD_2)\Th_2^*\1n\\
\noalign{\smallskip}\displaystyle \qq-(\Th_1^*)^T\big(D_1^TPD_2\1n+\1nR_{12}\big)
\Th_2^*\1n-\1n(\Th_2^*)^T\big(D_2^TPD_1\1n+\1nR_{21}\big)\Th_1^*\\
\noalign{\smallskip}\displaystyle \q\1n=PA\1n+\1nA^TP\1n+\1nC^TPC\1n+\1nQ\1n-\1n(\Th^*)^T(R\1n+\1nD^TPD)\Th^*\\
\noalign{\smallskip}\displaystyle \q\1n=PA\1n+\1nA^TP\1n+\1nC^TPC\1n+\1nQ\1n-
\1n\big(PB\1n+\1nC^TPD\1n+\1nS^T\big)(R\1n+\1nD^TPD)^\dag\big(B^TP\1n+\1nD^TPC\1n+\1nS\big).\ea\ee
Therefore, $P$ is a stabilizing solution of ARE (\ref{R-G}). \endpf

\ms

The following result, which is the main result of this paper, gives a characterization for closed-loop saddle points of Problem (LQG).

\ms

\bf Theorem 5.7. \sl Problem {\rm(LQG)} admits a closed-loop saddle
point $(\Th^*,u^*(\cd))\1n\in\dbR^{m\times n}\1n\times\1n
L^2_\dbF(\dbR^m)$ with $\Th^*\equiv((\Th_1^*)^T,(\Th_2^*)^T)^T$ and
$u^*(\cd)\equiv(u_1^*(\cd)^T,u_2^*(\cd)^T)^T$  if and only if the
following hold:

\ms

{\rm(i)} ARE $(\ref{R-G})$ admits a stabilizing solution $P$;

\ms

{\rm(ii)} The following BSDE:
\bel{eta-zeta}\ba{ll}
\noalign{\smallskip}\displaystyle  d\eta\1n=\1n-\Big\{\big[A^T\1n-\1n\sL(P)\sN(P)^\dag B^T\,\big]\eta\1n+\1n\big[C^T\1n-\1n\sL(P)\sN(P)^\dag D^T\,\big]\z\\
\noalign{\smallskip}\displaystyle \qq\q\ \,+\big[C^T\1n-\1n\sL(P)\sN(P)^\dag D^T\,\big]P\si\1n-\1n\sL(P)\sN(P)^\dag\rho\1n+\1nPb\1n+\1nq\Big\}dt\1n+\1n\z dW(t),
\q t\ges0,\ea\ee
admits an $L^2$-stable adapted solution $(\eta(\cd),\z(\cd))$ such
that
\bel{}B^T\eta(t)+D^T\z(t)+D^TP\si(t)+\rho(t)\in\cR\big(\sN(P)\big),\q \ae t\in[0,\infty),~\as\ee
In this case, the closed-loop saddle point $(\Th^*,u^*(\cd))$ admits
the following representation:
\bel{representation}\left\{\2n\ba{ll}
\noalign{\smallskip}\displaystyle \Th^*\1n=-\sN(P)^\dag\sL(P)^T\1n+\1n\big[I\1n-\sN(P)^\dag\sN(P)\big]\Pi,\\
\noalign{\smallskip}\displaystyle
u^*(\cd)\1n=-\sN(P)^\dag\big[B^T\eta(\cd)\1n+\1nD^T\z(\cd)\1n+\1nD^TP\si(\cd)\1n+\1n\rho(\cd)\big]\1n
+\1n\big[I\1n-\1n\sN(P)^\dag\sN(P)\big]\n(\cd),\ea\right.\ee
where $\Pi\in\dbR^{m\times n}$ is chosen such that
$\Th^*\in\sS[A,C;B,D]$, and $\n(\cd)\in L_\dbF^2(\dbR^m)$.

\ms

Further, the value function admits the following representation:
\bel{Value}\ba{ll}
\noalign{\smallskip}\displaystyle  V(x)=\lan Px,x\ran\1n+\dbE\,\Big\{2\lan\eta(0),x\ran\1n+\1n\int_0^\infty
\3n\big[\lan P\si,\si\ran\1n+2\lan\eta,b\ran\1n+2\lan\z,\si\ran\\
\noalign{\smallskip}\displaystyle \qq\qq\qq\q-\1n \lan(R\1n+\1n D^T\1nPD)^\dag(B^T\1n\eta\1n+\1nD^T\1n\z\1n+\1n D^T\1n P\si\1n+\1n\rho),
B^T\1n\eta\1n+\1nD^T\1n\z\1n+\1nD^T\1n P\si\1n+\1n\rho\ran\big]dt\Big\}.\ea\ee

\ms

\it Proof. \rm {\it Necessity.} Let
$(\Th^*,u^*(\cd))\1n\in\1n\dbR^{m\times n}\1n\times\1n
L^2_\dbF(\dbR^m)$ be a closed-loop saddle point of Problem (LQG)
with $\Th^*\equiv((\Th_1^*)^T,(\Th_2^*)^T)^T$ and
$u^*(\cd)\equiv(u_1^*(\cd)^T,u_2^*(\cd)^T)^T$. It follows from
Proposition 5.4 that $(\Th_1^*,0;\Th_2^*,0)$ is a closed-loop saddle
point of Problem $\hb{\rm{(LQG)}}^0$. By Proposition 5.6, ARE
(\ref{R-G}) admits a stabilizing solution $P$, and $\Th^*$ is given
by (\ref{K^*}).

\ms

To determine $u^*(\cd)$, let $(X^*(\cd), Y^*(\cd), Z^*(\cd))$ be the solution of (\ref{FBSDE-G}). Then
\bel{u-1}Ru^*+B^TY^*+D^TZ^*+(S+R\Th^*)X^*+\rho=0,\q \ae~\as\ee
and hence,
$$\ba{ll}
\noalign{\smallskip}\displaystyle  dY^*\1n=\1n-\big\{(A\1n+\1nB\Th^*)^TY^*\1n+(C\1n+\1nD\Th^*)^TZ^*\1n+\1n\wt QX^*\1n
+\1n\wt S^Tu^*\1n+\1n\wt q\,\big\}dt\1n+\1nZ^*dW(t)\\
\noalign{\smallskip}\displaystyle \qq=\1n-\big\{A^TY^*\1n+\1nC^TZ^*\1n+\1n(Q\1n+\1nS^T\Th^*)X^*\1n+\1nS^Tu^*\1n+\1nq\\
\noalign{\smallskip}\displaystyle \qq\qq\,+(\Th^*)^T\big[B^TY^*\1n+\1nD^TZ^*\1n+\1n(S+R\Th^*)X^*\1n+\1n
Ru^*\1n+\1n\rho\big]\big\}dt\1n+\1nZ^*dW(t)\\
\noalign{\smallskip}\displaystyle \qq=\1n-\big\{A^TY^*\1n+\1nC^TZ^*\1n+\1n(Q\1n+\1nS^T\Th^*)X^*\1n+\1nS^Tu^*\1n+\1nq\big\}dt\1n+\1nZ^*dW(t),\q t\ges0.\ea$$
Define
$$\left\{\2n\ba{ll}
\noalign{\smallskip}\displaystyle \eta(t)=Y^*(t)-PX^*(t),\\
\noalign{\smallskip}\displaystyle \z(t)=Z^*(t)-P(C+D\Th^*)X^*(t)-PDu^*(t)-P\si(t),\ea\right.\qq
t\ges0.$$
Noting $\sM(P)+\sL(P)\Th^*=0$, we have
$$\ba{ll}
\noalign{\smallskip}\displaystyle  \q\2n d\eta=dY^*-PdX^*\\
\noalign{\smallskip}\displaystyle \qq=-\big[A^TY^*+C^TZ^*+(Q+S^T\Th^*)X^*+S^Tu^*+q\big]dt+Z^*dW\\
\noalign{\smallskip}\displaystyle \qq\q-P\big[(A+B\Th^*)X^*+Bu^*+b\big]dt-P\big[(C+D\Th^*)X^*+Du^*+\si\big]dW\\
\noalign{\smallskip}\displaystyle \qq=-\Big\{A^T(\eta+PX^*)+C^T\big[\z+P(C+D\Th^*)X^*+PDu^*+P\si\big]\\
\noalign{\smallskip}\displaystyle \qq\qq\q+(Q+S^T\Th^*)X^*+S^Tu^*+q+P\big[(A+B\Th^*)X^*+Bu^*+b\big]\Big\}dt+\z dW\\
\noalign{\smallskip}\displaystyle \qq=-\Big\{
A^T\1n\eta+C^T\1n\z\1n+\1n\sM(P)X^*\1n+\1n\sL(P)\Th^*X^*\1n
+\1n\sL(P)u^*\1n+\1n C^T\1n P\si\1n+\1n Pb\1n+\1n q\Big\}dt+\z dW\\
\noalign{\smallskip}\displaystyle \qq=-\big[A^T\eta+C^T\z+\sL(P)u^*+C^TP\si+Pb+q\big]dt+\z dW.\ea$$
According to (\ref{u-1}), we have (noting $\sL(P)^T+\sN(P)\Th^*=0$)
$$\ba{ll}
\noalign{\smallskip}\displaystyle 0=B^TY^*+D^TZ^*+(S+R\Th^*)X^*+Ru^*+\rho\\
\noalign{\smallskip}\displaystyle \q\1n=B^T(\eta\1n+\1nPX^*)\1n+\1nD^T\big[\z\1n+\1nP(C\1n+\1nD\Th^*)X^*\1n+\1nPDu^*\1n+\1nP\si\big]\1n
+\1n(S\1n+\1nR\Th^*)X^*\1n+\1nRu^*\1n+\1n\rho\\
\noalign{\smallskip}\displaystyle \q\1n=\big[\sL(P)^T+\sN(P)\Th^*\big]X^*+B^T\eta+D^T\z+D^TP\si+\rho+\sN(P)u^*\\
\noalign{\smallskip}\displaystyle \q\1n=B^T\eta+D^T\z+D^TP\si+\rho+\sN(P)u^*.\ea$$
Hence,
$$B^T\eta+D^T\z+D^TP\si+\rho\in\cR\big(\sN(P)\big),\q\ae~\as$$
Since  $\sN(P)^\dag(B^T\eta\1n+\1nD^T\z\1n+\1nD^TP\si\1n+\1n\rho)\1n=\1n-\sN(P)^\dag\sN(P)u^*,$ and $\sN(P)^\dag\sN(P)$
is an orthogonal projection, we have
$$u^*=-\sN(P)^\dag(B^T\eta+D^T\z+D^TP\si+\rho)+\big[I-\sN(P)^\dag\sN(P)\big]\n$$
for some $\n(\cd)\in L_\dbF^2(\dbR^m)$. Consequently,
$$\ba{ll}
\noalign{\smallskip}\displaystyle  \sL(P)u^*=-\sL(P)\sN(P)^\dag(B^T\eta\1n+\1nD^T\z\1n+\1nD^TP\si\1n+\1n\rho)\1n+\1n\sL(P)\big[I\1n-\1n\sN(P)^\dag\sN(P)\big]\n\\
\noalign{\smallskip}\displaystyle \qq\qq=-\sL(P)\sN(P)^\dag(B^T\eta\1n+\1nD^T\z\1n+\1nD^TP\si\1n+\1n\rho).\ea$$
Then
$$\ba{ll}
\noalign{\smallskip}\displaystyle  A^T\eta+C^T\z+\sL(P)u^*+C^TP\si+Pb+q\\
\noalign{\smallskip}\displaystyle \,=A^T\eta+C^T\1n\z-\sL(P)\sN(P)^\dag(B^T\1n\eta+D^T\1n\z+
D^T\1n P\si+\rho)+C^T\1n P\si+Pb+q\\
\noalign{\smallskip}\displaystyle \,=\big[A^T-\sL(P)\sN(P)^\dag B^T\,\big]\eta+\big[C^T-\sL(P)\sN(P)^\dag D^T\,\big]\z\\
\noalign{\smallskip}\displaystyle \qq+\big[C^T-\sL(P)\sN(P)^\dag D^T\,\big]P\si-\sL(P)\sN(P)^\dag\rho+Pb+q.\ea$$
Therefore, $(\eta,\z)$ is an $L^2$-stable solution to (\ref{eta-zeta}).

\ms

{\it Sufficiency.} Let $(\Th^*,u^*(\cd))$ be given by
(\ref{representation}), where $\Pi\in\dbR^{m\times n}$ is chosen so
that $\Th^*\in\sS[A,C;B,D]$. Then
\bel{Ito-1}\sN(P)\Th^*\1n+\1n\sL(P)^T\1n=\1n0,\q
\sM(P)\1n+\1n\sL(P)\Th^*\1n+\1n(\Th^*)^T\sL(P)^T\1n+\1n(\Th^*)^T\sN(P)\Th^*\1n=\1n0,\ee
\bel{Ito-2}B^T\eta+D^T\z+D^TP\si+\rho=-\sN(P)u^*,\ee
and
\bel{Ito-3}\big[(\Th^*)^T\1n+\1n\sL(P)\sN(P)^\dag\big](B^T\eta\1n+\1nD^T\z\1n+\1nD^TP\si\1n+\1n\rho)
\1n=\1n-\Pi^T\big[I\1n-\1n\sN(P)\sN(P)^\dag\big]\sN(P)u^*\1n=\1n0.\ee
We take any
$u(\cd)\1n=\1n(u_1(\cd)^T\1n,u_2(\cd)^T)^T\1n\in\1n L^2_\dbF(\dbR^{m_1})\1n\times L^2_\dbF(\dbR^{m_2})$, and let
$X(\cd)\1n\equiv \1nX(\cd\,;x,u(\cd))$ be the solution of the following closed-loop system:
$$\left\{\2n\ba{ll}
\noalign{\smallskip}\displaystyle
dX(t)\1n=\1n\Big\{\1n\big[A\1n+\1nB\Th^*\big]X(t)\1n+\1nBu(t)\1n+\1nb(t)\1n\Big\}dt\1n
+\1n\Big\{\1n\big[C\1n+\1nD\Th^*\big]X(t)\1n+\1nDu(t)\1n+\1n\si(t)\1n\Big\}dW(t),\q t\ges 0,\\
\noalign{\smallskip}\displaystyle  X(0)\1n=\1nx.\ea\right.$$
Then
\bel{J-1}\ba{ll}
\noalign{\smallskip}\displaystyle  J(x;\Th^*X(\cd)\1n+\1nu(\cd))\1n=\1n\dbE\1n\int_0^\infty\1n\Big[\1n\lan\1n\begin{pmatrix}Q&S^T\\
                                                  S&R\end{pmatrix}
                                   \begin{pmatrix}X\\
                                                  \Th^*X\1n+\1nu\end{pmatrix},
                                   \begin{pmatrix}X\\
                                                  \Th^*X\1n+\1nu\end{pmatrix}\1n\ran
                                                  \1n+2\lan\1n\begin{pmatrix}q\\
                                                    \rho\end{pmatrix},
                                   \begin{pmatrix}X\\
                                                  \Th^*X\1n+\1nu\end{pmatrix}\1n\ran\1n\Big]dt\\
\noalign{\smallskip}\displaystyle \,=\dbE\int_0^\infty\1n\Big\{\lan\big[Q+S^T\Th^*+(\Th^*)^TS+(\Th^*)^TR\Th^*\big]X,X\ran
+2\lan(S+R\Th^*)X,u\ran\\
\noalign{\smallskip}\displaystyle \qq\qq\q+\lan Ru,u\ran+2\lan q+(\Th^*)^T\rho,X\ran+2\lan
\rho,u\ran\Big\}dt.\ea\ee
Applying It\^o's formula to $t\mapsto\lan PX(t),X(t)\ran$, one has (noting (\ref{Ito-1}))
\bel{J-2}\ba{ll}
\noalign{\smallskip}\displaystyle \lan
Px,x\ran\1n=\1n-\dbE\int_0^\infty\2n\Big\{\1n\lan\big[P(A\1n+\1nB\Th^*)\1n+\1n(A\1n+\1nB\Th^*)^TP\big]X,X\ran\1n
+\1n\lan P(C\1n+\1nD\Th^*)X,(C\1n+\1nD\Th^*)X\ran\\
\noalign{\smallskip}\displaystyle \qq\qq\qq+2\lan PX,Bu\1n+\1nb\ran\1n+2\lan P(C\1n+\1nD\Th^*)X,Du\1n+\1n\si\ran)\1n
+\1n\lan P(Du\1n+\1n\si),Du\1n+\1n\si\ran\1n\Big\}dt\\
\noalign{\smallskip}\displaystyle \,=\1n-\dbE\int_0^\infty\2n\Big\{\1n\lan\big[(PA\1n+\1nA^TP\1n+\1nC^TPC)\1n
+\1n(PB\1n+\1nC^TPD)\Th^*\1n+\1n(\Th^*)^T(B^TP\1n+\1nD^TPC)\\
\noalign{\smallskip}\displaystyle \qq\qq\q+(\Th^*)^TD^TPD\Th^*\big]X,X\ran\1n+2\lan(B^TP\1n+\1nD^TPC\1n+\1nD^TPD\Th^*)X,u\ran\1n\\
\noalign{\smallskip}\displaystyle \qq\qq\q+2\lan P(C\1n+\1nD\Th^*)X,\si\ran\1n+\1n\lan D^TPDu,u\ran\1n+2\lan D^TP\si,u\ran\1n
+2\lan PX,b\ran\1n+\1n\lan P\si,\si\ran \1n\Big\}dt\\
\noalign{\smallskip}\displaystyle \,=\1n-\dbE\int_0^\infty\2n\Big\{\1n\lan\big[\sM(P)\1n+\1n\sL(P)\Th^*\1n+\1n(\Th^*)^T\sL(P)^T\1n
+\1n(\Th^*)^T\sN(P)\Th^*\big]X,X\ran\\
\noalign{\smallskip}\displaystyle \qq\qq\q-\lan\big[Q\1n+\1nS^T\Th^*\1n+\1n(\Th^*)^TS\1n+\1n(\Th^*)^TR\Th^*\big]X,X\ran\\
\noalign{\smallskip}\displaystyle \qq\qq\q+2\lan\big[\sL(P)^T\1n+\1n\sN(P)\Th^*\1n-\1n(S\1n+\1nR\Th^*)\big]X,u\ran\1n\\
\noalign{\smallskip}\displaystyle \qq\qq\q+2\lan P(C\1n+\1nD\Th^*)X,\si\ran\1n+\1n\lan
D^TPDu,u\ran\1n
+2\lan D^TP\si,u\ran\1n+2\lan PX,b\ran\1n+\1n\lan P\si,\si\ran \1n\Big\}dt\\
\noalign{\smallskip}\displaystyle \,=\1n-\dbE\int_0^\infty\2n\Big[2\lan P(C\1n+\1nD\Th^*)X,\si\ran\1n
+\1n\lan D^TPDu,u\ran\1n+2\lan D^TP\si,u\ran\1n+2\lan PX,b\ran\1n+\1n\lan P\si,\si\ran \Big]dt\\
\noalign{\smallskip}\displaystyle \q\,+\dbE\int_0^\infty\1n\lan\big[Q\1n+\1nS^T\Th^*\1n+\1n(\Th^*)^TS\1n+\1n(\Th^*)^TR\Th^*\big]X,X\ran\1n
+2\lan(S\1n+\1nR\Th^*)X,u\ran dt.\ea\ee
Applying It\^o's formula to $t\mapsto\lan \eta(t),X(t)\ran$, one has (noting (\ref{Ito-3}))
\bel{J-3}\ba{ll}
\noalign{\smallskip}\displaystyle  \dbE\lan\1n \eta(0),x\1n\ran\1n=\1n\dbE\1n\int_0^\infty\2n\Big\{\1n\lan\1n\big[A^T\1n-\1n\sL(P)\sN(P)^\dag B^T\,\big]\eta\1n+\1n\big[C^T\1n-\1n\sL(P)\sN(P)^\dag D^T\,\big]\z\\
\noalign{\smallskip}\displaystyle \qq\qq\qq\qq\q+\big[C^T\1n-\1n\sL(P)\sN(P)^\dag D^T\,\big]P\si\1n-\1n\sL(P)\sN(P)^\dag\rho\1n+\1nPb\1n+\1nq,X\1n\ran\\
\noalign{\smallskip}\displaystyle \qq\qq\qq\qq\q-\1n\lan(A\1n+\1nB\Th^*)X\1n+\1nBu\1n+\1nb,\eta\1n\ran\1n-\1n\lan\1n\z,(C\1n+\1nD\Th^*)X\1n
+\1nDu\1n+\1n\si\1n\ran\1n\Big\}dt\\
\noalign{\smallskip}\displaystyle \qq\qq\q\1n=\1n\dbE\1n\int_0^\infty\2n\Big\{\1n-\1n\lan\1n\big[(\Th^*)^T\1n+\1n\sL(P)\sN(P)^\dag\big]B^T\eta\1n
+\1n\big[(\Th^*)^T\1n+\1n\sL(P)\sN(P)^\dag \big]D^T\z,X\1n\ran\\
\noalign{\smallskip}\displaystyle \qq\qq\qq\qq\q\1n-\1n\lan\1n\big[(\Th^*)^T\1n+\1n\sL(P)\sN(P)^\dag\big]D^TP\si,X\1n\ran\1n
+\1n\lan\1n P(C\1n+\1nD\Th^*)X,\si\1n\ran\\
\noalign{\smallskip}\displaystyle \qq\qq\qq\qq\q\1n-\1n\lan\1n\sL(P)\sN(P)^\dag\rho,X\1n\ran\1n+\1n\lan\1n Pb\1n
+\1nq,X\1n\ran\1n-\1n\lan\1n Bu\1n+\1nb,\eta\1n\ran\1n-\1n\lan\1n\z,Du\1n+\1n\si\1n\ran\1n\Big\}dt\\
\noalign{\smallskip}\displaystyle \qq\qq\q\1n=\1n\dbE\1n\int_0^\infty\2n\Big\{\1n-\1n\lan\1n\big[(\Th^*)^T\1n+\1n\sL(P)\sN(P)^\dag\big](B^T\eta\1n
+\1nD^T\z\1n+\1nD^TP\si\1n+\1n\rho),X\1n\ran\\
\noalign{\smallskip}\displaystyle \qq\qq\qq\qq\q\1n+\1n\lan\1n P(C\1n+\1nD\Th^*)X,\si\1n\ran\1n+\1n\lan(\Th^*)^T\rho\1n+\1nPb\1n
+\1nq,X\1n\ran\1n-\1n\lan\1n Bu\1n+\1nb,\eta\1n\ran\1n-\1n\lan\1n\z,Du\1n+\1n\si\1n\ran\1n\Big\}dt\\
\noalign{\smallskip}\displaystyle \qq\qq\q\1n=\1n\dbE\1n\int_0^\infty\2n\Big\{\1n\lan\1n
P(C\1n+\1nD\Th^*)X,\si\1n\ran\1n+\1n\lan(\Th^*)^T\rho\1n+\1nPb\1n+\1nq,X\1n\ran\1n-\1n\lan\1n
Bu\1n+\1nb,\eta\1n\ran\1n-\1n\lan\1n\z,Du\1n+\1n\si\1n\ran\1n\Big\}dt.\ea\ee
Combining (\ref{J-1})--(\ref{J-3}) and noting (\ref{Ito-2}), we have
$$\ba{ll}
\noalign{\smallskip}\displaystyle  J(x;\Th^*X(\cd)\1n+\1nu(\cd))\1n-\1n\lan Px,x\ran\1n-2\,\dbE\lan \eta(0),x\ran\\
\noalign{\smallskip}\displaystyle \,=\dbE\1n\int_0^\infty\1n\Big\{\1n\lan\sN(P)u,u\ran\1n
+2\lan B^T\eta\1n+\1nD^T\z\1n+\1nD^TP\si\1n+\1n\rho,u\ran\1n+2\lan b,\eta\ran\1n+2\lan \z,\si\ran\1n+\1n\lan P\si,\si\ran\1n\Big\}dt\\
\noalign{\smallskip}\displaystyle \,=\dbE\1n\int_0^\infty\1n\Big\{\1n\lan\sN(P)u,u\ran\1n-2\lan \sN(P)u^*,u\ran\1n+2\lan b,\eta\ran\1n+2\lan\z,\si\ran\1n+\1n\lan P\si,\si\ran\1n\Big\}dt\\
\noalign{\smallskip}\displaystyle \,=\dbE\1n\int_0^\infty\1n\Big\{\1n\lan\sN(P)(u\1n-\1nu^*),u\1n-\1nu^*\ran\1n-\1n\lan \sN(P)u^*,u^*\ran\1n
+2\lan b,\eta\ran\1n+2\lan\z,\si\ran\1n+\1n\lan P\si,\si\ran\1n\Big\} dt.\ea$$
Consequently,
$$\ba{ll}
\noalign{\smallskip}\displaystyle  J(x;\Th_1^*X(\cd)+u_1(\cd),\Th_2^*X(\cd)+u_2^*(\cd))-J(x;\Th^*X^*(\cd)+u^*(\cd))\\
\noalign{\smallskip}\displaystyle \,=\dbE\int_0^\infty\1n\lan(R_{11}+D_1^TPD_1)(u_1-u_1^*),u_1-u_1^*\ran dt\ges0\ea$$
since $R_{11}+D_1^TPD_1\ges0$. Similarly,
$$\ba{ll}
\noalign{\smallskip}\displaystyle  J(x;\Th_1^*X(\cd)+u_1^*(\cd),\Th_2^*X(\cd)+u_2(\cd))-J(x;\Th^*X^*(\cd)+u^*(\cd))\\
\noalign{\smallskip}\displaystyle \,=\dbE\int_0^\infty\1n\lan(R_{22}+D_2^TPD_2)(u_2-u_2^*),u_2-u_2^*\ran dt\les0\ea$$
since $R_{22}+D_2^TPD_2\les0$. Therefore, $(\Th^*,u^*(\cd))$ is a
closed-loop saddle point of Problem (LQG). Finally, noting
(\ref{Ito-2}), we have
$$\ba{ll}
\noalign{\smallskip}\displaystyle  \lan\sN(P)u^*,u^*\ran=\lan\sN(P)\sN(P)^\dag\sN(P)u^*,u^*\ran=\lan\sN(P)^\dag\sN(P)u^*,\sN(P)u^*\ran\\
\noalign{\smallskip}\displaystyle \qq\qq\qq\q\2n=\lan(R\1n+\1n D^T\1nPD)^\dag(B^T\1n\eta\1n+D^T\1n\z\1n+\1n D^T\1n P\si\1n+\1n\rho),
B^T\1n\eta\1n+D^T\1n\z\1n+\1nD^T\1n P\si\1n+\1n\rho\ran,\ea$$
and hence,
$$\ba{ll}
\noalign{\smallskip}\displaystyle  V(x)\1n=\1nJ(x;\Th^*X(\cd)\1n+\1nu^*(\cd))\\
\noalign{\smallskip}\displaystyle \qq\ =\1n\lan Px,x\ran\1n+2\,\dbE\1n\lan \eta(0),x\ran\1n+\dbE\int_0^\infty\2n\Big\{\2n-\1n\lan \sN(P)u^*,u^*\ran\1n
+2\lan b,\eta\ran\1n+2\lan\z,\si\ran\1n+\1n\lan P\si,\si\ran\1n\Big\} dt\\
\noalign{\smallskip}\displaystyle \qq\ =\1n\lan Px,x\ran\1n+\dbE\,\Big\{2\lan\eta(0),x\ran\1n+\1n\int_0^\infty
\3n\big[\lan P\si,\si\ran\1n+2\lan\eta,b\ran+2\lan\z,\si\ran\\
\noalign{\smallskip}\displaystyle \qq\qq\qq\q-\1n \lan(R\1n+\1n D^T\1nPD)^\dag(B^T\1n\eta\1n+D^T\1n\z\1n+\1n D^T\1n P\si\1n+\1n\rho),
B^T\1n\eta\1n+D^T\1n\z\1n+\1nD^T\1n P\si\1n+\1n\rho\ran\big]dt\Big\}.\ea$$
This completes the proof. \endpf

\ms

Note that the above result is reduced to that for Problem (LQ) if
$m_2=0$. It is not hard for us to state such a result and we omit
the details here.

\section{Examples}

In this section we present two examples illustrating how the ``stabilizing solution" of AREs plays an important role in the
study of closed-loop saddle points. For simplicity, we only consider one player (optimal control) case.

\ms

\bf Example 6.1. \rm Consider the following state equation
$$\left\{\2n\ba{ll}
\noalign{\smallskip}\displaystyle  dX(t)=-\big[2X(t)+u(t)\big]dt+\big[2X(t)+u(t)\big]dW(t),\q t\ges0,\\
\noalign{\smallskip}\displaystyle  X(0)=x,\ea\right.$$
with the cost functional
$$J(x;u(\cd))=\dbE\int^\infty_0\2n\[\,2\,|X(t)|^2-{1\over2}|u(t)|^2\,\]dt.$$
By Lemma 2.3, part (iv), the system $[-2,2;-1,1]$ is stabilizable, and
$\Th\in\sS[-2,2;-1,1]$ if and only if
$$2(-2-\Th)+(2+\Th)^2<0\q (\hb{\,i.e.,\,}-2<\Th<0).$$
The corresponding ARE reads
$$P^2-2P+1=0.$$
Thus, $P=1$ and
$$\big[I-\sN(P)^\dag\sN(P)\big]\Pi-\sN(P)^\dag\sL(P)^T\equiv-2,\qq\forall\,\Pi\in\dbR.$$
Hence, by Theorem 5.7, the above problem does not admit any closed-loop optimal control. From this example, we see that
ARE (\ref{R-G}) may only admit non-stabilizing solutions.

\ms

\bf Example 6.2. \rm Consider the following state equation
$$\left\{\2n\ba{ll}
\noalign{\smallskip}\displaystyle  dX(t)=-\[\,{1\over4}X(t)+2u(t)\]dt+\big[X(t)+u(t)\big]dW(t),\q t\ges0,\\
\noalign{\smallskip}\displaystyle  X(0)=x,\ea\right.$$
with the cost functional
$$J(x;u(\cd))=\dbE\int^\infty_0\2n\[\,{1\over2}|X(t)|^2-2X(t)u(t)+|u(t)|^2\,\]dt.$$
By Lemma 2.3, part (iv), $\Th\in\sS[-{1\over4},1;-2,1]$ if and only if
$$2\(\1n-{1\over4}-2\Th\)+(1+\Th)^2<0\q \(\hb{\,i.e., }1-{\sqrt{2}\over2}<\Th<1+{\sqrt{2}\over2}\;\).$$
The corresponding ARE reads
$$(P+1)^2=0,$$
which admits a unique stabilizing solution $P=-1$. Noting $\sN(P)=0$, by Theorem 5.7, we see that
$$(\Pi,\n(\cd));\qq \Pi\in\(1-{\sqrt{2}\over2},1+{\sqrt{2}\over2}\;\),\ \n(\cd)\in L_\dbF^2(\dbR)$$
are all the closed-loop optimal controls of the above problem. However,
$$-\sN(P)^\dag\sL(P)^T=0\not\in\sS\[-{1\over4},1;-2,1\].$$
Also, from this example, we see that even if $-\sN(P)^\dag\sL(P)^T$ is not a
stabilizer of the system, Problem (LQ) may still admit closed-loop optimal controls.


\end{document}